\documentclass[12pt]{amsart}
\usepackage{amssymb,amsmath, amsthm,latexsym,verbatim}
\usepackage{a4wide}
\usepackage[hypertex]{hyperref}
\usepackage{enumitem}
\usepackage{colortbl}

\newcommand{\cal}[1]{\mathcal{#1}}
\theoremstyle{plain}
\newtheorem{lem}{Lemma}[section]
\newtheorem{theo}[lem]{Theorem}
\newtheorem{prop}[lem]{Proposition}
\newtheorem{corollary}[lem]{Corollary}
\newtheorem{remark}[lem]{Remark}
\newtheorem{definition}{Definition}
\newtheorem{conjecture}{Conjecture}

\parskip=\medskipamount
\font\k=cmr7
\font\rm=cmr12

  \newcommand {\free}{\mbox{\k free}}
  \newcommand {\tors}{\mbox{\k tors}}
  \newcommand {\topo}{\mbox{\k top}}
  
  \newcommand {\cu}{\mbox{\k cus}}
  \newcommand {\di}{\mbox{\k disc}}
  \newcommand {\sd}{\mbox{\k sd}}
  \newcommand {\sph}{\mbox{\k sph}}
  \newcommand {\fin}{\mbox{\k fin}}
  \newcommand {\pl}{\text{pl}}

  \newcommand {\reg}{\mbox{\k reg}}
  \newcommand {\res}{\mbox{\k res}}
  \newcommand {\spec}{\mbox{\k spec}}
  \newcommand {\geo}{\mbox{\k geo}}
  
  \newcommand {\temp}{\text{temp}}

  \newcommand {\C}{{\mathbb C}}
  \newcommand {\bH}{{\mathbb H}}
  \newcommand {\N}{{\mathbb N}}
  \newcommand {\R}{{\mathbb R}}
  \newcommand {\Z}{{\mathbb Z}}
  \newcommand {\Q}{{\mathbb Q}}
  \newcommand {\A}{{\mathbb A}}

  \newcommand {\af}{{\mathfrak a}}
  \newcommand {\gf}{{\mathfrak g}}
  
  \newcommand {\kf}{{\mathfrak k}}

  \newcommand {\ho}{{\mathfrak o}}
  \newcommand {\pf}{{\mathfrak p}}

\renewcommand {\H}{{\mathcal H}}
  \newcommand {\cM}{{\mathcal M}}
  
  \newcommand {\cO}{{\mathcal O}}
  \newcommand {\G}{{\bf G}}
  \newcommand {\K}{{\bf K}}
  \newcommand {\bP}{{\bf P}}
  \newcommand {\bQ}{{\bf Q}}
  \newcommand {\bU}{{\bf U}}
  \newcommand {\M}{{\bf M}}
  \newcommand {\bL}{{\bf L}}
  \newcommand {\bR}{{\bf R}}

  \newcommand {\E}{{\mathcal E}}
  \newcommand {\cD}{{\mathcal D}}
 \newcommand {\cP}{{\mathcal P}}
 
 \newcommand {\cL}{{\mathcal L}}

\newcommand {\bs}{\backslash}

\newcommand  {\T}{{\bf S}}
\renewcommand{\Im}{\operatorname{Im}}
\renewcommand{\Re}{\operatorname{Re}}
\newcommand{\Tr}{\operatorname{Tr}}

\newcommand{\End}{\operatorname{End}}

\newcommand{\inter}{\operatorname{int}}
\newcommand{\injrad}{\operatorname{injrad}}
\newcommand{\tr}{\operatorname{tr}}
\newcommand{\Id}{\operatorname{Id}}
\newcommand{\Hom}{\operatorname{Hom}}

\newcommand{\Ind}{\operatorname{Ind}}

\newcommand{\rk}{\operatorname{rank}}
\newcommand{\vol}{\operatorname{vol}}
\newcommand{\Area}{\operatorname{Area}}

\newcommand{\SL}{\operatorname{SL}}
\newcommand{\GL}{\operatorname{GL}}
\newcommand{\SO}{\operatorname{SO}}

\newcommand{\Spin}{\operatorname{Spin}}

\newcommand{\inj}{\operatorname{inj}}

\renewcommand{\det}{\operatorname{det}}

\newcommand{\Sym}{\operatorname{Sym}}

\newcommand{\norm}[1]{\lVert#1\rVert}
\newcommand{\AAA}{A}
\newcommand{\levis}{{\mathcal L}}
\newcommand{\aaa}{\mathfrak{a}}
\newcommand{\Ht}{H}
\newcommand{\sprod}[2]{\left\langle#1,#2\right\rangle}
\newcommand{\abs}[1]{\lvert#1\rvert}
\newcommand{\PPP}{\mathcal{P}}
\newcommand{\FFF}{{\mathcal F}}
\newcommand{\rts}{\Sigma}
\newcommand{\disc}{\operatorname{disc}}
\newcommand{\srts}{\Delta}
\newcommand{\modulus}{\delta}
\newcommand{\AF}{{\mathcal A}}
\newcommand{\zzz}{\mathfrak{z}}
\newcommand{\iii}{{\mathrm i}}
\newcommand{\LieG}{\mathfrak{g}}
\newcommand{\unip}{\operatorname{unip}}
\newcommand{\bases}{\mathfrak{B}}
\newcommand{\bss}{\underline{\beta}}
\newcommand{\dtup}{\mathcal{X}}
\newcommand{\card}[1]{\lvert#1\rvert}

\newcommand{\level}{\operatorname{level}}

\newcommand{\param}{\Lambda}
\newcommand{\TWN}{(TWN)}
\newcommand{\BD}{(BD)}

\newcommand{\types}{\mathcal{F}}
\newcommand{\eps}{\epsilon}

\newcommand{\data}{\mathcal{D}}

\begin{document}
\title[]
{Asymptotics of automorphic spectra and the trace formula}

\author{Werner M\"uller}
\address{Universit\"at Bonn\\
Mathematisches Institut\\
Beringstrasse 1\\
D -- 53115 Bonn, Germany}
\email{mueller@math.uni-bonn.de}

\begin{abstract}
This paper is a survey article on the limiting behavior  of the discrete 
spectrum of the right regular
representation in $L^2(\Gamma\bs G)$ for a lattice $\Gamma$ in a reductive
group $G$ over a number field. We discuss various aspects of the Weyl law,
the limit multiplicity problem and the analytic torsion. 
\end{abstract}

\maketitle
\setcounter{tocdepth}{1}
\tableofcontents
\section{Introduction}
Let $G$ be a connected, linear, semisimple algebraic group over $\Q$. Let 
$\Pi(G(\R))$ denote
the set of all equivalence classes of irreducible unitary representations of
$G(\R)$, equipped with the Fell topology \cite{Di}. We fix a Haar measure 
on $G(\R)$. Let $\Gamma\subset G(\R)$ be a lattice, i.e., a discrete subgroup 
such 
that $\vol(\Gamma\bs G(\R))<\infty$. Let $R_\Gamma$ be the right regular 
representation of $G$ on $L^2(\Gamma\bs G)$. Let $L^2_{\di}(\Gamma\bs G)$ be
the span of all irreducible subrepresentations of $R_\Gamma$ and denote by
$R_{\Gamma,\di}$ the restriction of $R_\Gamma$ to $L^2_{\di}(\Gamma\bs G)$. Then
$R_{\Gamma,\di}$ decomposes discretely as
\begin{equation}
R_{\Gamma,\di}\cong\widehat\bigoplus_{\pi\in\Pi(G)}m_\Gamma(\pi)\pi,
\end{equation}
where
\[
m_\Gamma(\pi)=\dim\Hom_G(\pi,R_\Gamma)=\dim\Hom_G(\pi,R_{\Gamma,\di})
\]
is the multiplicity with which $\pi$ occurs in $R_\Gamma$. The multiplicities are
known to be finite under a weak reduction-theoretic assumption on $(G,\Gamma)$
\cite{OW}, which is satisfied if $G$ has no compact factors or if $\Gamma$ is
arithmetic. The study of the multiplicities $m_\Gamma(\pi)$ is one of the main
concerns in the theory of automorphic forms. Apart from special cases like
discrete series representations, one cannot hope in general to describe the
multiplicity function on $\Pi(G)$ explicitly. A more feasible and 
interesting problem is the study of
the asymptotic behavior of the multiplicities with respect to the growth of
various parameters such as the level of congruence subgroups or the 
infinitesimal character of $\pi$. 
This is closely related to the study of families of automorphic forms
 (see \cite{SST}).

The first problem in this context is the Weyl law. Let $K$ be a maximal
compact subgroup of $G$. Fix an irreducible representation $\sigma$ of $K$.
Let $\Pi(G;\sigma)$ be the subspace of all $\pi\in\Pi(G)$ such that $[\pi|_K
\colon \sigma]>0$. Especially, if $\sigma_0$ is the trivial representation,
then $\Pi(G;\sigma_0)$ is the spherical dual $\Pi(G(\R))_{\sph}$. Given $\pi\in
\Pi(G(\R))$, denote by $\lambda_\pi=\pi(\Omega)$ the Casimir eigenvalue of 
$\pi$.
For $\lambda\ge 0$ let the counting function be defined by
\begin{equation}\label{counting1}
N_\Gamma(\lambda;\sigma)=
\sum_{\substack{\pi\in\Pi(G;\sigma)\\|\lambda_\pi|\le\lambda}}m_\Gamma(\pi).
\end{equation}
Then the problem is to determine the behavior of the counting function as
$\lambda\to\infty$.

Another basic problem is the limit multiplicity problem, which is the study of 
the asymptotic behavior of the multiplicities if $\vol(\Gamma\bs G(\R))\to 
\infty$. For $G=\GL_n$ this corresponds
to the study of harmonic families of cuspidal automorphic representations
of $\GL_n(\A)$ (see \cite{SST}). 
More precisely, for a given lattice $\Gamma$ define
the discrete spectral measure $\mu_\Gamma$ on $\Pi(G)$, associated to 
$\Gamma$, by
\begin{equation}
\mu_\Gamma=\frac{1}{\vol(\Gamma\bs G(\R))}\sum_{\pi\in\Pi(G(\R))}
m_\Gamma(\pi)\delta_\pi,
\end{equation}
where $\delta_\pi$ is the Dirac measure at $\pi$. Then the limit multiplicity
problem is concerned with the study of the asymptotic behavior of 
$\mu_\Gamma$ as $\vol(\Gamma\bs G(\R))\to \infty$. For appropriate sequences of
lattices $(\Gamma_n)$ one expects that the measures $\mu_{\Gamma_n}$ converge to
the Plancherel measure $\mu_{\pl}$ on $\Pi(G(\R))$.

There are closely related problems in topology and spectral theory. One of
them concerns Betti numbers.
Let $K$ be a maximal compact subgroup of $G$ and put $X=G/K$. Let
$\Gamma$ be a uniform lattice in $G$ and let $(\Gamma_n)$ be a tower of
normal subgroups of $\Gamma$. Put $M=\Gamma\bs X$ and $M_n=
\Gamma_n\bs X$, $n\in\N$.  Then $M_n\to M$ is a sequence of finite
normal coverings of $M$.  For any topological space $Y$ let $b_k(Y)$
denote the $k$-th Betti number of  $Y$. Then
\begin{equation}\label{betti}
\lim_{n\to\infty}\frac{b_k(M_n)}{\vol(M_n)}=b^{(2)}_k(X),
\end{equation}
where $b^{(2)}_k(X)$ is the $k$-th $L^2$-Betti number of $X$. This was proved by 
L\"uck \cite{Lu1} in the more general context of CW-complexes. In the case of 
locally 
symmetric spaces, it follows from the results about limit multiplicities. 
Again, it was extended by Abert et al \cite{AB1} to much more general 
sequences of uniform lattices.

A more sophisticated spectral invariant is the Ray-Singer 
analytic torsion $T_X(\rho)$ (see \cite{RS}). It depends on a finite 
dimensional  representation $\rho$ of $\Gamma$ and is defined in terms of the
spectra of the Laplace operators $\Delta_p(\rho)$ on $p$-forms with coefficients
in the flat bundle associated with $\rho$. Of particular interest are 
representations of $\Gamma$ which arise as the 
restriction of a representation of $G$. For appropriate representations, called
strongly acyclic, Bergeron and Venkatesch \cite{BV} studied the asymptotic 
behavior of $\log T_{X_n}(\rho)$ as $n\to\infty$. One of their main results is
\begin{equation}\label{lim-tor}
\lim_{n\to\infty}\frac{\log T_{X_n}(\rho)}{\vol(X_n)}=\log T^{(2)}_X(\rho),
\end{equation}
where $T^{(2)}_X(\rho)$ is the $L^2$-torsion \cite{Lo}, \cite{MV}. 
Using the equality of analytic torsion and Reidemeister torsion \cite{Ch},
\cite{Mu1}, \eqref{lim-tor} implies results about the growth of 
 the torsion subgroup in the integer homology of arithmetic groups. Let 
${\bf G}$ be a semisimple algebraic group over $\Q$, $G={\bf G}(\R)$ and
$\Gamma\subset {\bf G}(\Q)$ a co-compact, arithmetic subgroup. As shown in
\cite{BV}, there are strongly acyclic representations $\rho$ of $G$ on
a finite dimensional vector space $V$ such that $V$ contains a 
$\Gamma$-invariant lattice $M$. Let $\cM$ be the  local system of 
free $\Z$-modules over $X$, attached to $M$. Then the cohomology $H_*(X,\cM)$ 
of $X$ with coefficients in $\cM$ is a finite abelian group. Denote by
$|H_*(X,\cM)|$ its order. Assume that $d=\dim(X)$ is odd. Then by \cite{BV} one has
\[
\lim_{n\to\infty}\sum_{p=1}^d(-1)^{p+\frac{d-1}{2}}\frac{\log |H_p(X_n,\cM)|}
{[\Gamma\colon\Gamma_n]}
=c_{M,G}\vol(X),
\]
where $c_{M,G}$ is a constant that depends only on $G$ and $M$. Moreover, if
$\delta(G):=\rk G-\rk K=1$, then $c_{M,G}>0$. It is conjectured that the limit
\begin{equation}\label{homol-tor}
\lim_{n\to\infty}\frac{\log |H_j(X_n,\cM)|}{[\Gamma\colon\Gamma_n]}
\end{equation}
always exists and is equal to zero, unless $\delta(G)=1$ and $j=(d-1)/2$. 
In the latter case it is equal to $c_{M,G}$ times $\vol(X)$. The conjecture is
known to be true for $G=\SL_2(\C)$. 

An important problem is to extend these results to the non-compact case.

\section{The Arthur trace formula}
\setcounter{equation}{0}

The trace formula is one of the main technical tools to study the kind of
spectral problems mentioned in the introduction. For $\R$-rank one groups
the Selberg trace formula is available \cite{Wa1}. In the higher rank case
the Selberg trace formula is replaced by the Arthur trace formula. 

In this section we recall Arthur's trace formula, and in particular
the refinement of the spectral expansion obtained in \cite{FLM1}.

\subsection{Notation} \label{subsecnotation}

We will mostly use the notation of \cite{FLM1}.
Let $\G$ be a reductive group defined over $\Q$ and let $\A$ be the ring of 
adeles of $\Q$.
We fix a maximal compact subgroup $\K=\prod_v \K_v = \K_\infty\cdot 
\K_{\fin}$ of $\G(\A)=\G(\R)\cdot \G(\A_{\fin})$.

Let $\gf$ and $\kf$ denote the Lie algebras of $\G(\R)$ and $\K_\infty$,
respectively. Let $\theta$ be the Cartan involution of $\G(\R)$ with respect to
$\K_\infty$. It induces a Cartan decomposition $\mathfrak{g}= 
\mathfrak{p} \oplus \mathfrak{k}$. 
We fix an invariant bi-linear form $B$ on $\mathfrak{g}$ which is positive 
definite on $\mathfrak{p}$ and negative definite on $\mathfrak{k}$.
This choice defines a Casimir operator $\Omega$ on $\G(\R)$,
and we denote the Casimir eigenvalue of any $\pi \in \Pi (\G(\R))$ by 
$\lambda_\pi$. Similarly, we obtain
a Casimir operator $\Omega_{\K_\infty}$ on $\K_\infty$ and write $\lambda_\tau$ for 
the Casimir eigenvalue of a
representation $\tau \in \Pi (\K_\infty)$ (cf.~\cite[\S 2.3]{BG}).
The form $B$ induces a Euclidean scalar product $(X,Y) = - B (X,\theta(Y))$ on 
$\mathfrak{g}$ and all its subspaces.
For $\tau \in \Pi (\K_\infty)$ we define $\norm{\tau}$ as in 
\cite[\S 2.2]{CD}. 

We fix a maximal $\Q$-split torus $\T_0$ of $\G$ and let
$\M_0$ be its centralizer, which is a minimal Levi subgroup defined over $\Q$.
We assume that the maximal compact subgroup $\K \subset \G(\A)$ is admissible 
with respect to $\M_0$ \cite[\S 1]{Ar5}.
Denote by $\AAA_0$ the identity component of $\T_0(\R)$, which is viewed as a 
subgroup of $\T_0(\A)$. We write $\levis$ for the (finite) set of Levi 
subgroups containing $\M_0$, i.e., the set of centralizers of subtori of $\T_0$.
Let $W_0=N_{\G(\Q)}(\T_0)/\M_0$ be the Weyl group of $(\G,\T_0)$,
where $N_{\G(\Q)}(H)$ is the normalizer of $H$ in $\G(\Q)$.
For any $s\in W_0$ we choose a representative $w_s\in \G(\Q)$.
Note that $W_0$ acts on $\levis$ by $s\M=w_s\M w_s^{-1}$.

Let now $\M\in\levis$. We write $\T_M$ for the split part of the identity 
component of the center of $\M$.
Set $\AAA_M=\AAA_0\cap \T_M(\R)$ and $W(\M)=N_{\G(\Q)}(\M)/\M$, which can be 
identified with a subgroup of $W_0$.
Denote by $\aaa_M^*$ the $\R$-vector space spanned by the lattice $X^*(\M)$ of 
$\Q$-rational characters of $\M$ and let $\aaa_{M,\C}^*=\aaa_M^*\otimes_{\R}\C$ 
be its complexification.
We write $\aaa_M$ for the dual space of $\aaa_M^*$, which is spanned by the 
co-characters of $\T_M$.
Let $\Ht_M:\M(\A)\rightarrow\aaa_M$ be the homomorphism given by
\[
e^{\sprod{\chi}{\Ht_M(m)}}=\abs{\chi (m)}_\A = \prod_v\abs{\chi(m_v)}_v
\]
for any $\chi\in X^*(\M)$ and denote by $\M(\A)^1 \subset \M (\A)$ the kernel 
of $\Ht_M$. Let $\levis(\M)$ be the set of Levi subgroups containing $\M$ and
$\PPP(\M)$ the set of parabolic subgroups of $\G$ with Levi part $\M$.
We also write $\FFF(\M)=\FFF^G(\M)=\coprod_{\bL\in\levis(\M)}\PPP(\bL)$ for the 
(finite) set of parabolic subgroups of $\G$ containing $\M$.
Note that $W(\M)$ acts on $\PPP(\M)$ and $\FFF(\M)$ by $s\bP=w_s\bP w_s^{-1}$.
Denote by $\rts_M$ the set of reduced roots of $\T_M$ on the Lie algebra of 
$\G$.
For any $\alpha\in\rts_M$ we denote by $\alpha^\vee\in\aaa_M$ the corresponding 
co-root. Let $L^2_{\disc}(\AAA_M\M(\Q)\bs \M(\A))$ be the discrete part of 
$L^2(\AAA_M\M(\Q)\bs \M(\A))$, i.e., the
closure of the sum of all irreducible subrepresentations of the regular 
representation of $\M(\A)$.
We denote by $\Pi_{\disc}(\M(\A))$ the countable set of equivalence classes of 
irreducible unitary
representations of $\M(\A)$ which occur in the decomposition of 
$L^2_{\disc}(\AAA_M\M(\Q)\bs \M(\A))$
into irreducible representations.

For any $\bL\in\levis(\M)$ we identify $\aaa_L^*$ with a subspace of $\aaa_M^*$.
We denote by $\aaa_M^L$ the annihilator of $\aaa_L^*$ in $\aaa_M$. We set
\[
\levis_1(\M)=\{\bL\in\levis(\M):\dim\aaa_M^L=1\}
\]
and
\[
\FFF_1(\M)=\bigcup_{\bL\in\levis_1(M)}\PPP(\bL).
\]
Note that
the restriction of the scalar product
$(\cdot,\cdot)$ on $\mathfrak{g}$ defined above
gives
$\aaa_{M_0}$ the structure of a Euclidean space.
In particular, this fixes Haar measures on the spaces $\aaa^L_M$ and their duals
$(\aaa^L_M)^*$. We follow Arthur in the corresponding normalization of Haar 
measures on the groups $\M(\A)$
(\cite[\S 1]{Ar1}).

\subsection{Intertwining operators}

The main ingredient of the spectral side of the Arthur trace formula are 
logarithmic derivatives of intertwining operators. We shall now describe the
structure of the intertwining operators.

Let $\bP\in\PPP(\M)$. We write $\aaa_P=\aaa_M$.
Let $\bU_P$ be the unipotent radical of $\bP$ and $\M_P$ the unique 
$\bL\in\levis(\M)$
(in fact the unique $\bL\in\levis(\M_0)$) such that $\bP\in\PPP(\bL)$.
Denote by $\rts_P\subset\aaa_P^*$ the set of reduced roots of $\T_M$ on the 
Lie algebra $\mathfrak{u}_P$ of $\bU_P$.
Let $\srts_P$ be the subset of simple roots of $\bP$, which is a basis for 
$(\aaa_P^G)^*$.
Write $\aaa_{P,+}^*$ for the closure of the Weyl chamber of $\bP$, i.e.
\[
\aaa_{P,+}^*=\{\lambda\in\aaa_M^*:\sprod{\lambda}{\alpha^\vee}\ge0
\text{ for all }\alpha\in\rts_P\}
=\{\lambda\in\aaa_M^*:\sprod{\lambda}{\alpha^\vee}\ge0\text{ for all }
\alpha\in\srts_P\}.
\]
Denote by $\modulus_P$ the modulus function of $\bP(\A)$.
Let $\bar\AF_2(P)$ be the Hilbert space completion of
\[
\{\phi\in C^\infty(\M(\Q)\bU_P(\A)\bs \G(\A)):\modulus_P^{-\frac12}\phi(\cdot x)\in
L^2_{\disc}(A_M\M(\Q)\bs \M(\A)),\ \forall x\in \G(\A)\}
\]
with respect to the inner product
\[
(\phi_1,\phi_2)=\int_{\AAA_M\M(\Q)\bU_P(\A)\bs \G(\A)}\phi_1(g)\overline{\phi_2(g)}\ dg.
\]


Let $\alpha\in\rts_M$.
We say that two parabolic subgroups $\bP,\bQ\in\PPP(\M)$ are \emph{adjacent} along $\alpha$,
and write $\bP|^\alpha \bQ$, if $\rts_P\cap-\rts_Q=\{\alpha\}$.
Alternatively, $\bP$ and $\bQ$ are adjacent if the closure $\overline{\bP\bQ}$ 
of $\bP\bQ$ belongs to $\FFF_1(\M)$.
Any ${\bf R}\in\FFF_1(\M)$ is of the form $\overline{\bP\bQ}$ for a unique
unordered pair $\{\bP,\bQ \}$ of parabolic subgroups in $\PPP(\M)$, namely 
$\bP$ and $\bQ$ are the maximal parabolic subgroups
of ${\bf R}$, and $\bP|^\alpha \bQ$ with $\alpha^\vee\in\rts_P^\vee \cap\aaa^R_M$.
Switching the order of $\bP$ and $\bQ$ changes $\alpha$ to $-\alpha$.

For any $\bP\in\PPP(\M)$ let $\Ht_P\colon\G(\A)\rightarrow\aaa_P$ be the 
extension of $\Ht_M$ to a left $\bU_P(\A)$-
and right $\K$-invariant map.
Denote by $\AF^2(P)$ the dense subspace of $\bar\AF^2(P)$ consisting of its 
$\K$- and $\zzz$-finite vectors,
where $\zzz$ is the center of the universal enveloping algebra of $\mathfrak{g} \otimes \C$.
That is, $\AF^2(P)$ is the space of automorphic forms $\phi$ on 
${\bf U}_P(\A){\bf M}(F)\bs \G(\A)$ such that
$\modulus_P^{-\frac12}\phi(\cdot k)$ is a square-integrable automorphic form on
$\AAA_MM(F)\bs M(\A)$ for all $k\in\K$.
Let $\rho(P,\lambda)$, $\lambda\in\aaa_{M,\C}^*$, be the induced
representation of $G(\A)$ on $\bar\AF^2(P)$ given by
\[
(\rho(P,\lambda,y)\phi)(x)=\phi(xy)e^{\sprod{\lambda}{\Ht_P(xy)-\Ht_P(x)}}.
\]
It is isomorphic to $\Ind_{\bP(\A)}^{\G(\A)}\left(L^2_{\disc}(\AAA_M \M(\Q)\bs \M(\A))\otimes e^{\sprod{\lambda}{\Ht_M(\cdot)}}\right)$.

For $\bP,\bQ\in\PPP(\M)$ let
\[
M_{Q|P}(\lambda):\AF^2(P)\to\AF^2(Q),\quad\lambda\in\aaa_{M,\C}^*,
\]
be the standard \emph{intertwining operator} \cite[\S 1]{Ar3}, which is the 
meromorphic
continuation in $\lambda$ of the integral
\[
[M_{Q|P}(\lambda)\phi](x)=\int_{\bU_Q(\A)\cap \bU_P(\A)\bs \bU_Q(\A)}\phi(nx)
e^{\sprod{\lambda}{\Ht_P(nx)-\Ht_Q(x)}}\ dn, \quad \phi\in\AF^2(P), \ x\in \G(\A).
\]
These operators satisfy the following properties.
\begin{enumerate}
\item $M_{P|P}(\lambda)\equiv\Id$ for all $\bP\in\PPP(\M)$ and 
$\lambda\in\aaa_{M,\C}^*$.
\item For any $\bP,\bQ,\bR\in\PPP(\M)$ we have 
$M_{R|P}(\lambda)=M_{R|Q}(\lambda)\circ M_{Q|P}(\lambda)$ for all 
$\lambda\in\aaa_{M,\C}^*$.
In particular, $M_{Q|P}(\lambda)^{-1}=M_{P|Q}(\lambda)$.
\item $M_{Q|P}(\lambda)^*=M_{P|Q}(-\overline{\lambda})$ for any 
$\bP,\bQ\in\PPP(\M)$
and $\lambda\in\aaa_{M,\C}^*$. In particular, $M_{Q|P}(\lambda)$
is unitary for $\lambda\in\iii\aaa_M^*$.
\item If $\bP|^\alpha \bQ$ then $M_{Q|P}(\lambda)$ depends only on 
$\sprod{\lambda}{\alpha^\vee}$.
\end{enumerate}
Given $\pi\in\Pi_{\di}(\M(\A))$, let $\AF^2_\pi(P)$ be the space of all 
$\phi\in\AF^2(P)$ for which the function 
$x\in \M(\A)\mapsto \delta_P^{-\frac{1}{2}}\phi(xg)$,
$g\in \G(\A)$, belongs to the $\pi$-isotypic subspace 
$L^2(A_M\M(\Q)\bs \M(\A))$.
For any $\bP\in\PPP(\M)$ we have a canonical isomorphism of 
$\G(\A_f)\times(\LieG_{\C},K_\infty)$-modules
\[
j_P:\Hom(\pi,L^2(\AAA_M\M(\Q)\bs \M(\A)))\otimes
\Ind_{\bP(\A)}^{\G(\A)}(\pi)\rightarrow\AF^2_\pi(P).
\]
If we fix a unitary structure on $\pi$ and endow 
$\Hom(\pi,L^2(\AAA_M\M(\Q)\bs \M(\A)))$ with the inner product 
$(A,B)=B^\ast A$
(which is a scalar operator on the space of $\pi$), the isomorphism $j_P$ 
becomes an isometry. 

Suppose that $\bP|^\alpha \bQ$.
The operator $M_{Q|P}(\pi,s):=M_{Q|P}(s\varpi)|_{\AF^2_\pi(P)}$, where $\varpi\in
\af°\ast_M$ is such that $\langle\varpi,\alpha^\vee\rangle=1$, admits a normalization by a global factor
$n_\alpha(\pi,s)$ which is a meromorphic function in $s$. We may write
\begin{equation} \label{normalization}
M_{Q|P}(\pi,s)\circ j_P=n_\alpha(\pi,s)\cdot j_Q\circ(\Id\otimes R_{Q|P}(\pi,s))
\end{equation}
where $R_{Q|P}(\pi,s)=\otimes_v R_{Q|P}(\pi_v,s)$ is the product
of the locally defined normalized intertwining operators and 
$\pi=\otimes_v\pi_v$
\cite[\S 6]{Ar3}, (cf.~\cite[(2.17)]{Mu6}). In many cases, the 
normalizing factors can be expressed in terms automorphic $L$-functions 
\cite{Sha1}, \cite{Sha2}. For example, let $\G=\GL(n)$. Then
the global normalizing factors $n_\alpha$ can be expressed in terms
of Rankin-Selberg $L$-functions and the
known properties of these functions, which are collected and analyzed in 
\cite[\S\S 4,5]{Mu5}.
Write $\M \simeq \prod_{i=1}^r \GL (n_i)$, where the root $\alpha$ is trivial on 
$\prod_{i \ge 3} \GL (n_i)$,
and let $\pi \simeq \otimes \pi_i$ with representations $\pi_i \in 
\Pi_{\disc}(\GL(n_i,\A))$.
Let $L(s,\pi_1\times\tilde\pi_2)$ be the completed Rankin-Selberg $L$-function 
associated to $\pi_1$ and $\pi_2$. It satisfies the functional equation
\begin{equation}\label{functequ}
L(s,\pi_1\times\tilde\pi_2)=\eps(\frac12,\pi_1\times\tilde\pi_2)
N(\pi_1\times\tilde\pi_2)^{\frac12-s}L(1-s,\tilde\pi_1\times\pi_2)
\end{equation}
where $\abs{\eps(\frac12,\pi_1\times\tilde\pi_2)}=1$ and $N(\pi_1\times\tilde
\pi_2)\in\N$ is the conductor. Then we have
\begin{equation}\label{rankin-selb}
n_\alpha (\pi,s) = \frac{L(s,\pi_1\times\tilde\pi_2)}{\eps(\frac12,\pi_1\times
\tilde\pi_2)N(\pi_1\times\tilde\pi_2)^{\frac12-s}L(s+1,\pi_1\times\tilde\pi_2)}.
\end{equation}

\subsection{The trace formula}
\setcounter{equation}{0}

Arthur's trace formula gives two alternative expressions for a distribution 
$J$ on $G (\A)^1$. Note that this
distribution depends on the choice of $M_0$ and $\K$.
For $h \in C^\infty_c (G (\A)^1)$, Arthur defines $J(h)$ as the
value at the point $T=T_0$ specified in \cite[Lemma 1.1]{Ar5} of a 
polynomial $J^T (h)$ on $\aaa_{M_0}$ of degree at most $d_0 = \dim \aaa^G_{M_0}$.
Here, the polynomial $J^T(h)$ depends in addition on the choice of a 
parabolic subgroup $P_0 \in \PPP (M_0)$.
Consider the equivalence relation on $G(\Q)$ defined by $\gamma \sim \gamma'$ 
whenever the semisimple parts of $\gamma$ and
$\gamma'$ are $G(\Q)$-conjugate.
Let $\mathcal{O}$ be the set of the resulting equivalence classes (which are 
in bijection with conjugacy classes of semisimple elements).
The coarse geometric expansion \cite{Ar1} is
\begin{equation} \label{geometricside}
J^T (h) = \sum_{\mathfrak{o} \in \mathcal{O}} J^T_\mathfrak{o} (h),
\end{equation}
where the summands $J^T_\mathfrak{o} (h)$ are again polynomials in $T$ of degree 
at most $d_0$.
Write $J_\mathfrak{o} (h) = J^{T_0}_\mathfrak{o} (h)$, which depends only on $M_0$ 
and $\K$.
Then
$J_\mathfrak{o} (h) = 0$ if the support of $h$ is disjoint from all conjugacy 
classes of $G(\A)$
intersecting
$\mathfrak{o}$ (cf.~\cite[Theorem 8.1]{Ar6}). By [ibid., Lemma 9.1] 
(together with the
descent formula of \cite[\S 2]{Ar5}), for each compact set
$\Omega \subset G (\A)^1$ there exists a finite subset $\mathcal{O} (\Omega) 
\subset \mathcal{O}$ such that
for $h$ supported in $\Omega$ only the terms with $\mathfrak{o} \in 
\mathcal{O} (\Omega)$
contribute to \eqref{geometricside}. In particular, the sum is always finite.
When $\mathfrak{o}$ consists of the unipotent elements of $G(\Q)$, we write 
$J^T_{\unip} (h)$ for $J^T_\mathfrak{o} (h)$.

We now turn to the spectral side. Let $L \supset M$ be Levi subgroups in 
$\levis$, $P \in \PPP (M)$, and
let $m=\dim\aaa_L^G$ be the co-rank of $L$ in $G$.
Denote by $\bases_{P,L}$ the set of $m$-tuples $\bss=(\beta_1^\vee,\dots,
\beta_m^\vee)$
of elements of $\rts_P^\vee$ whose projections to $\aaa_L$ form a basis for 
$\aaa_L^G$.
For any $\bss=(\beta_1^\vee,\dots,\beta_m^\vee)\in\bases_{P,L}$ let
$\vol(\bss)$ be the co-volume in $\aaa_L^G$ of the lattice spanned by $\bss$ 
and let
\begin{align*}
\Xi_L(\bss)&=\{(\bQ_1,\dots,\bQ_m)\in\FFF_1(M)^m: \ \ \beta_i^{\vee}\in\aaa_M^{Q_i}, 
\, i = 1, \dots, m\}\\&=
\{(\overline{P_1P_1'},\dots,\overline{P_mP_m'}): \ \ P_i|^{\beta_i}P_i', \, 
i = 1, \dots, m\}.
\end{align*}

For any smooth function $f$ on $\aaa_M^*$ and $\mu\in\aaa_M^*$ denote by 
$D_\mu f$ the directional derivative of $f$ along $\mu\in\aaa_M^*$.
For a pair $P_1|^\alpha P_2$ of adjacent parabolic subgroups in $\PPP(M)$ write
\[
\delta_{P_1|P_2}(\lambda)=M_{P_2|P_1}(\lambda)D_\varpi M_{P_1|P_2}(\lambda):
\AF^2(P_2)\rightarrow\AF^2(P_2),
\]
where $\varpi\in\aaa_M^*$ is such that $\sprod{\varpi}{\alpha^\vee}=1$.
\footnote{Note that this definition differs slightly from the definition of
$\delta_{P_1|P_2}$ in \cite{FL1}.} Equivalently, writing
$M_{P_1|P_2}(\lambda)=\Phi(\sprod{\lambda}{\alpha^\vee})$ for a
meromorphic function $\Phi$ of a single complex variable, we have
\[
\delta_{P_1|P_2}(\lambda)=\Phi(\sprod{\lambda}{\alpha^\vee})^{-1}
\Phi'(\sprod{\lambda}{\alpha^\vee}).
\]
For any $m$-tuple $\dtup=(Q_1,\dots,Q_m)\in\Xi_L(\bss)$
with $Q_i=\overline{P_iP_i'}$, $P_i|^{\beta_i}P_i'$, denote by 
$\Delta_{\dtup}(P,\lambda)$
the expression
\[
\frac{\vol(\bss)}{m!}M_{P_1'|P}(\lambda)^{-1}\delta_{P_1|P_1'}(\lambda)M_{P_1'|P_2'}(\lambda) \cdots
\delta_{P_{m-1}|P_{m-1}'}(\lambda)M_{P_{m-1}'|P_m'}(\lambda)\delta_{P_m|P_m'}(\lambda)M_{P_m'|P}(\lambda).
\]
In \cite[pp. 179-180]{FLM1} we define a (purely combinatorial) map $\dtup_L: \bases_{P,L} \to \FFF_1 (M)^m$ with the property that
$\dtup_L(\bss) \in \Xi_L (\bss)$ for all $\bss \in \bases_{P,L}$.\footnote{The map $\dtup_L$ depends in fact on the additional choice of
a vector $\underline{\mu} \in (\aaa^*_M)^m$ which does not lie in an explicit finite
set of hyperplanes. For our purposes, the precise definition of $\dtup_L$ is immaterial.}

For any $s\in W(M)$ let $L_s$ be the smallest Levi subgroup in $\levis(M)$
containing $w_s$. We recall that $\aaa_{L_s}=\{H\in\aaa_M\mid sH=H\}$.
Set
\[
\iota_s=\abs{\det(s-1)_{\aaa^{L_s}_M}}^{-1}.
\]
For $P\in\FFF(M_0)$ and $s\in W(M_P)$ let
$M(P,s):\AF^2(P)\to\AF^2(P)$ be as in \cite[p.~1309]{Ar3}.
$M(P,s)$ is a unitary operator which commutes with the operators $\rho(P,\lambda,h)$ for $\lambda\in\iii\aaa_{L_s}^*$.
Finally, we can state the refined spectral expansion.

\begin{theo}[\cite{FLM1}] \label{thm-specexpand}
For any $h\in C_c^\infty(G(\A)^1)$ the spectral side of Arthur's trace formula is given by
\begin{equation}\label{specside1}
J (h) = \sum_{[M]} J_{\spec,M} (h),
\end{equation}
$M$ ranging over the conjugacy classes of Levi subgroups of $G$ (represented by members of $\mathcal{L}$),
where
\begin{equation}\label{specside2}
J_{\spec,M} (h) =
\frac1{\card{W(M)}}\sum_{s\in W(M)}\iota_s
\sum_{\bss\in\bases_{P,L_s}}\int_{\iii(\aaa^G_{L_s})^*}
\tr(\Delta_{\dtup_{L_s}(\bss)}(P,\lambda)M(P,s)\rho(P,\lambda,h))\ d\lambda
\end{equation}
with $P \in \PPP(M)$ arbitrary.
The operators are of trace class and the integrals are absolutely convergent.
\end{theo}
Note that the term corresponding to $M=G$ is $J_{\spec,G} (h) = \tr R_{\disc}(h)$. 
Next assume that $M$ is the Levi subgroup of a maximal parabolic subgroup $P$.
Furthermore, let $L=M$. Let $\bar P$ be the opposite parabolic subgroup to $P$. 
Then up to a constant, the contribution to the spectral side is given by
\[
\sum_{\pi\in\Pi_{\di}(M(\A)^1)}\int_{i\af^\ast}\tr(M_{\bar P|P}(\pi,\lambda)^{-1}
\frac{d}{dz}M_{\bar P|P}(\pi,\lambda)M(P,s)\rho(P,\pi,\lambda,h))\;d\lambda.
\]
Now assume that $G=\SL(2,\R)$ and $K=\SO(2)$. Let $h\in\SL(2,\R))$ be 
bi-K-invariant. Let $C(s)$ be the scattering matrix.

\section{The Weyl law}
\setcounter{equation}{0}

The Weyl law is concerned with the study of the asymptotic behavior of the
counting function \eqref{counting1} as $\lambda\to\infty$. This is the first
problem which needs to be solved in order to be able to pursue a deeper study
of the cuspidal automorphic spectrum. For example, the study of statistical 
properties of the automorphic spectrum requires first of all to know that the
spectrum is infinite and has the right asymptotic properties. This, in 
particular, concerns the study of families of automorphic forms 
(see \cite{SST}).

The investigation of the asymptotic behavior of the counting function 
\eqref{counting1} is
closely related to the study of the counting function of the
eigenvalues of the Laplace operator on a compact Riemannian manifold \cite{DG}. 
Let $\widetilde X=G/K$. It can be  equipped with a $G$-invariant metric which
is unique up to scaling. Let 
$X=\Gamma\bs\widetilde X$. Assume that $\Gamma$  is torsion free. Then $X$ is a 
complete Riemannian manifold of finite volume. Let $\sigma\in\widehat K$ and 
let $\widetilde 
E_\sigma\to \widetilde X$ be the homogeneous vector bundle associated to 
$\sigma$, which is equipped with the invariant Hermitian metric induced by
$\sigma$. Let $E_\sigma=\Gamma\bs\widetilde E_\sigma$ be the corresponding 
locally homogeneous vector bundle over $X$. Let $\nabla^\sigma$ be the 
connection in $E_\sigma$ induced by the canonical connection in 
$\widetilde E_\sigma$. Let $\Delta_\sigma=(\nabla^\sigma)^\ast\nabla^\sigma$
be the Bochner-Laplace operator, acting in $C^\infty(X,E_\sigma)$.
It is an elliptic, second order, formally self-adjoint differential operator 
of Laplace type, i.e., its principal symbol is given by 
$\|\xi\|_x^2\Id_{E_{\sigma,x}}$. 
The Bochner-Laplace operator is related to the Casimir operator 
$R_\Gamma(\Omega)$ by
\begin{equation}\label{bochlapl}
\Delta_\sigma=-R_\Gamma(\Omega)+\lambda_\sigma\Id,
\end{equation}
where $\lambda_\sigma$ is the Casimir eigenvalue of $\sigma$. Assume that $X$
is compact. Then $\Delta_\sigma$ has a pure discrete spectrum consisting of a 
sequence of eigenvalues $0\le\lambda_1\le\lambda_2\le\cdots\to\infty$ of finite 
multiplicities. Let 
\[
N_{\Gamma}(\lambda;\sigma)=\#\{j\colon\lambda_j\le\lambda\}
\] 
be the counting function of the eigenvalues, where eigenvalues are counted  
with their multiplicity. By
\eqref{bochlapl} the counting function \eqref{counting1} has the same asymptotic
behavior as $N_{\Gamma}(\lambda;\sigma)$. The Weyl law for 
$N_{\Gamma}(\lambda;\sigma)$ can be
established by standard methods. For example, for a weak version, which means 
with no estimation of the remainder term, on can use the asymptotic expansion 
of the trace of the heat operator $e^{-t\Delta_\sigma}$. 
Thus if $\Gamma$ is co-compact we get from these general methods
the following formula for the asymptotic behavior of the counting function. Let
$d=\dim X$. As $\lambda\to\infty$ we have
\begin{equation}\label{weyl2}
N_{\Gamma}(\lambda;\sigma)=\frac{\dim(\sigma)\vol(\Gamma\bs G/K)}{(4\pi)^{d/2}
{\bf \Gamma}(d/2+1)}\lambda^{d/2}+o(\lambda^{d/2}),
\end{equation}
where ${\bf \Gamma}(s)$ denotes the Gamma function.

If $\Gamma$ is not co-compact, then $\Delta_\sigma$ has a nonempty continuous
spectrum which consists of a half-line $[c,\infty)$ for some $c\ge 0$. 
This makes it much more difficult to study the discrete spectrum of this
operator, because almost all eigenvalues, if they exist, will be embedded
into the continuous spectrum. It is well know from mathematical physics that
embedded eigenvalues are unstable under perturbations. One of the basic tools
to study the cuspidal automorphic spectrum is the trace formula.

\subsection{Hyperbolic surfaces}

In the non-compact case, a general Weyl law was first derived by Selberg for 
a hyperbolic surface $X=\Gamma\bs\bH$ of finite area, where 
$\bH=\SL(2,\R)/\SO(2)$ is the upper half-plane. We briefly recall the method
which is based on the trace formula. It illustrates the basic idea which is
also used in the higher rank case.

Let $\Delta=d^\ast d$ be the Laplace operator with respect to the hyperbolic 
metric. Then $\Delta$, regarded as operator in $L^2(X)$ with domain 
$C^\infty(X)$, is essentially self-adjoint. The spectrum of $\Delta$ is the union
of a pure point spectrum and the absolutely continuous spectrum. The pure point
spectrum consists of a sequence of eigenvalues
\[
0=\lambda_0<\lambda_1\le \lambda_2\le\cdots
\]
of finite multiplicities.  
If $X$ is noncompact then, in general, we only know that $\lambda_0$ exists.
We slightly change the definition of the counting function by
\[
N_\Gamma(\lambda):=\#\{j\colon\sqrt{\lambda_j}\le\lambda\}.
\]
The new terms in the trace formula, which are due to
the non-compactness of $\Gamma\bs\bH$ arise from the parabolic conjugacy
classes in $\Gamma$ and the Eisenstein series. Let us recall the definition of
Eisenstein series.  Let $a_1,...,a_m\in\R\cup\{\infty\}$ be
representatives of the $\Gamma$-conjugacy classes of parabolic fixed points of
$\Gamma$. The $a_i$'s are called {\it cusps}. For each $a_i$ let $\Gamma_{a_i}$
be the stabilizer of $a_i$ in $\Gamma$. Choose 
$\sigma_i\in\SL(2,\R)$ such that  
\[\sigma_i(\infty)=a_i,\quad \sigma_i^{-1}\Gamma_{a_i}\sigma_i=\left\{
\begin{pmatrix}1&n\\0&1\end{pmatrix}\colon n\in\Z\right\}.\]
Then the Eisenstein series $E_i(z,s)$ associated to the cusp $a_i$ is defined as
\begin{equation}
E_i(z,s)=\sum_{\gamma\in\Gamma_{a_i}\bs \Gamma}\Im(\sigma_i^{-1}\gamma z)^s,\quad 
\Re(s)>1.
\end{equation}
The series converges absolutely and uniformly on compact subsets of the 
half-plane $\Re(s)>1$ and it satisfies the following properties.
\begin{enumerate}
\item[1)] $E_i(\gamma z,s)=E_i(z,s)$ for all $\gamma\in\Gamma$.
\item[2)] As a function of $s$, $E_i(z,s)$ admits a meromorphic continuation 
to $\C$ which is regular on the line $\Re(s)=1/2$.
\item[3)] $E_i(z,s)$ is a smooth function of $z$ and satisfies
$\Delta_zE_i(z,s)=s(1-s) E_i(z,s).$
\end{enumerate}
The contribution of the Eisenstein series to the Selberg trace formula 
is given by their zeroth Fourier coefficients of the Fourier expansion in the 
cusps. The 
zeroth Fourier coefficient of the Eisenstein series $E_k(z,s)$ at the cusp 
$a_l$ is given by
\[\int_0^1 E_k(\sigma_l(x+iy),s)\;dx=\delta_{kl}y^s+C_{kl}(s)y^{1-s},\]
where $\delta_{kl}$ is Kronecker's delta function and $C_{kl}(s)$ is a 
meromorphic function of $s\in\C$. Put
\[
C(s):=\left(C_{kl}(s)\right)_{k,l=1}^m.
\]
This is the so called {\it scattering matrix}. Let $g\in C^\infty_c(\R)$ and
let $h=\hat g$ be the Fourier transform of $g$. Let $\phi(s):=\det C(s)$.
Denote by $\{\gamma\}$ the hyperbolic $\Gamma$-conjugacy classes. For every
hyperbolic element $\gamma$, denote by $\gamma_0$ the primitive hyperbolic 
element such that $\gamma=\gamma_0^k$ for some $k\in\N$. Every nontrivial
hyperbolic conjugacy class $\{\gamma\}$ corresponds to a unique closed 
geodesic $c_\gamma$. Let $l(\gamma)$ denote its length. Write the eigenvalues as
\[
\lambda_j=\frac{1}{4}+r_j^2,\quad r_j\in i\R\cup(1/2,1].
\]
Then the trace formula is the following identity.
\begin{equation}\label{traceform1}
\begin{split}
\sum_{j} h(r_j) - \frac{1}{4\pi}& \int^\infty_{-\infty}h(r)\frac{\phi^\prime}
{\phi}(1/2+ir)\;dr  +\frac{1}{4}\phi(1/2)h(0)\\
&=\frac{\Area(\Gamma\bs \bH)}{4\pi}\int_\R h(r)r\tanh(\pi r)\;dr
+\sum_{\{\gamma\}} 
\frac{l(\gamma_0)}{2\sinh \left( \frac{l(\gamma)}{2}\right)} g(l(\gamma))\\
& \quad -\frac{m}{2\pi}\int^\infty_{-\infty} h(r)
\frac{{\bf \Gamma}^\prime}{\bf \Gamma}(1+ir)dr +\frac{m}{4} h(0)-m\ln 2\; g(0).
\end{split}
\end{equation}
The left hand side is the spectral side, which contains all terms associated 
with the spectrum and the right hand side is the geometric side. 
The trace formula holds for every discrete subgroup $\Gamma\subset \SL(2,\R)$ 
with co-finite area. In analogy to the counting function of the eigenvalues we
introduce the winding number 
\begin{equation}\label{winding1}
M_\Gamma(\lambda)=-\frac{1}{4\pi}\int_{-\lambda}^\lambda\frac{\phi^\prime}
{\phi}(1/2+ir)\;dr,
\end{equation}
which measures the continuous spectrum. Using the cut-off Laplacian of 
Lax-Phillips \cite{CV} one can deduce the following elementary bounds 
\begin{equation}\label{bounds1}
N_\Gamma(\lambda)\ll\lambda^2,\quad M_\Gamma(\lambda)\ll\lambda^2,\quad 
\lambda\ge 1.
\end{equation}
These bounds imply that the the trace formula \eqref{traceform1} holds for a 
larger class of functions. In particular, it can be applied to the heat 
kernel $k_t$.  Its spherical Fourier transform 
equals $h_t(r)=e^{-t(1/4+r^2)}$, $t>0$.
If we insert $h_t$ into the trace formula we get the following asymptotic 
expansion as $t\to 0$.
\begin{equation}\label{asympexp1}
\begin{split}
\sum_{j}e^{-t\lambda_j}-\frac{1}{4\pi}\int_\R &e^{-t(1/4+r^2)}
\frac{\phi^\prime}{\phi}(1/2+ir)\;dr\\
&=\frac{\Area(\Gamma\bs\bH)}{4\pi t}+\frac{a\log t}{\sqrt{t}}+
\frac{b}{\sqrt{t}}+O(1)
\end{split}
\end{equation} 
for certain constants $a,b\in\R$. Using \cite[(8.8), (8.9)]{Se1} it 
follows that the winding number $M_\Gamma(\lambda)$ is monotonically increasing 
for $\lambda\gg0$. Therefore we can apply a Tauberian theorem to  
\eqref{asympexp1}
and we get the following  Weyl law, established by Selberg \cite{Se1}. 
As $\lambda\to\infty$ we have
\begin{equation}\label{weyllaw3}
N_\Gamma(\lambda)+M_\Gamma(\lambda)\sim \frac{\Area(\Gamma\bs\bH)}{4\pi}
\lambda^2.
\end{equation}
In general, we cannot estimate separately the counting function and the winding
number. For congruence subgroups, however, the entries of
the scattering matrix can be expressed in terms of well-known analytic
functions. For $\Gamma(N)$ the determinant of the scattering 
matrix $\phi(s)$ has been computed by Huxley \cite{Hu}. It has the 
form
\begin{equation}\label{4.5}
\phi(s)=(-1)^lA^{1-2s}\left(\frac{{\bf \Gamma}(1-s)}{{\bf \Gamma}(s)}\right)^k
\prod_\chi\frac{L(2-2s,\bar\chi)}{L(2s,\chi)},
\end{equation}
where $k,l\in\Z$, $A>0$, the product runs
over Dirichlet characters $\chi$ to some modulus dividing $N$ and $L(s,\chi)$
is the Dirichlet $L$-function with character $\chi$. Especially for $\Gamma(1)$
we have
\begin{equation}\label{4.5a}
\phi(s)=\sqrt{\pi}\frac{{\bf \Gamma}(s-1/2)\zeta(2s-1)}{{\bf \Gamma}(s)
\zeta(2s)},
\end{equation}
where $\zeta(s)$ denotes the Riemann zeta function.

Using Stirling's 
approximation formula to estimate the logarithmic derivative of the Gamma 
function and standard estimations for the logarithmic derivative
of Dirichlet $L$-functions on the line $\Re(s)=1$ 
\cite[Chapt V, Theormem 7.1]{Pr}, 
we get
\begin{equation}\label{4.6}
\frac{\phi'}{\phi}(1/2+ir)=O(\log(4+|r|)),\quad |r|\to\infty.
\end{equation}
This implies that
\begin{equation}\label{4.7}
M_{\Gamma(N)}(\lambda)\ll\lambda\log\lambda.
\end{equation}
Together with \eqref{weyllaw3} we obtain Weyl's law for the point spectrum of
the Laplacian on $X(N)=\Gamma(N)\bs\bH$:
\begin{equation}\label{4.8}
N_{\Gamma(N)}(\lambda) \sim
\frac{\Area(X(N))}{4\pi}\lambda^2,\quad\lambda\to\infty,
\end{equation}
which is due to Selberg \cite[p.668]{Se1}.
A similar formula holds for other congruence groups such as $\Gamma_0(N)$.
In particular,  (\ref{4.8})  implies
that for congruence groups there exist infinitely many linearly 
independent Maass cusp forms. 

By a more sophisticated use of the Selberg trace
formula one can estimate the remainder term (see \cite{Mu3}). For congruence
subgroups one gets
\begin{theo}
For every $N\in\N$ we have
\begin{equation}
N_{\Gamma(N)}(\lambda)=\frac{\Area(X(N))}{4\pi}\lambda^2+O(\lambda\log \lambda)
\end{equation}
as $\lambda\to\infty$.
\end{theo}
A finite area hyperbolic surface for which the Weyl law holds is called by
Sarnak {\it essentially cuspidal}. Now it is strongly believed that essential
cuspidality is limited to special arithmetic surfaces. This is based on work by 
Phillips and Sarnak who studied the behavior of the discrete spectrum when 
$\Gamma$ is deformed in the corresponding Teichm\"uller space. We refer to 
\cite{Sa1} for a detailed discussion of their method. This led Phillips and
Sarnak to the following conjectures.

\begin{conjecture}

\smallskip
1) The generic $\Gamma$ in a given Teichm\"uller space of finite area
hyperbolic surfaces is not essentially cuspidal.

\smallskip
\noindent
2) Except for the Teichm\"uller space of the once punctured torus, the generic
$\Gamma$ has only a finite number of discrete eigenvalues.
\end{conjecture}

\subsection{Higher rank}
We turn now to the general case. We assume that $G=\G(\R)$, where $\G$ is a 
connected semisimple algebraic group over $\Q$. Let $X=\Gamma\bs\widetilde X=
\Gamma\bs G/K$  and $E_\sigma\to X$ be as above. Let 
$\Delta_\sigma\colon C^\infty(X,E_\sigma)\to 
C^\infty(X,E_\sigma)$ be the Bochner-Laplace operator. As operator in 
$L^2(X,E_\sigma)$ it is essentially self-adjoint. 
Let $L^2_{\di}(X,E_\sigma)$ be the subspace of $L^2(X,E_\sigma)$ which is the
closure of the span of all $L^2$-eigensections of $\Delta_\sigma$. Recall that
a cusp form for $\Gamma$ is a smooth $K$-finite function $\phi\colon \Gamma\bs G
\to\C$ which is a joint eigenfunction of the center of the universal enveloping
algebra ${\cal Z}(\gf_\C)$ and which satisfies
\[
\int_{\Gamma\cap N_P\bs N_P}\phi(nx)\;dn=0
\]
for all unipotent radicals $N_P$ of proper rational parabolic subgroups $P$ of
$G$, i.e., $P=\bP(\R)$, where $\bP$ is a rational parabolic subgroup of $\G$.
Put
\[
L^2_{\cu}(X,E_\sigma):=(L^2_{\cu}(\Gamma\bs G)\otimes V_\sigma)^K.
\]
Then $L^2_{\cu}(X,E_\sigma)$ is contained in $L^2_{\di}(X,E_\sigma)$. 
The orthogonal complement $L^2_{\res}(X,E_\sigma)$  of $L^2_{\cu}(X,E_\sigma)$ in
$L^2_{\di}(X,E_\sigma)$ is called the {\it residual subspace}. By Langland's 
theory of Eisenstein series it follows that $L^2_{\res}(X,E_\sigma)$ is spanned 
by iterated residues of cuspidal Eisenstein series.
By definition we have an orthogonal decomposition
\[
L^2_{\di}(X,E_\sigma)=L^2_{\cu}(X,E_\sigma)\oplus L^2_{\res}(X,E_\sigma).
\]

Let $N_{\Gamma}^{\di}(\lambda;\sigma)$, $N_{\Gamma}^{\cu}(\lambda;\sigma)$, and
$N_{\Gamma}^{\res}(\lambda;\sigma)$ be the counting function of the eigenvalues
with eigensections belonging to the corresponding subspace. 
The following results about the growth of
the counting functions hold for any lattice $\Gamma$ in a real
semisimple Lie group. Let $d=\dim X$. Donnelly \cite{Do} has proved the
following bound for the cuspidal spectrum
\begin{equation}\label{limsup}
\limsup_{\lambda\to\infty}\frac{N_\Gamma^{\cu}(\lambda,\sigma)}{\lambda^{d/2}}\le
\frac{\dim(\sigma)\vol(X)}{(4\pi)^{d/2}
{\bf \Gamma}\left(\frac{d}{2}+1\right)}.
\end{equation}
For  the full discrete spectrum, we have at least an upper bound for the
growth of the counting function. The main result of \cite{Mu2} states that 
\begin{equation}\label{upperbd1}
N_\Gamma^{\di}(\lambda,\sigma)\ll(1+\lambda^{2d}).
\end{equation}
This result implies that invariant integral operators are of trace class on
the discrete subspace which is the starting point for the trace formula.
 The proof of \eqref{upperbd1} relies on the description
of the residual subspace in terms of iterated residues of Eisenstein series.

Let $N_\Gamma^{\cu}(\lambda)$ be the counting function with respect to the 
trivial representation $\sigma_0$ of $K$, i.e., the counting function of
the cuspidal spectrum of the Laplacian on functions. Then Sarnak \cite{Sa2} 
conjectured that if $\rk(G/K)>1$, Weyl's law holds for $N_\Gamma^{\cu}(\lambda)$,
which means that equality holds in \eqref{limsup}. Furthermore, one expects 
that 
the growth of the residual spectrum is of lower order than the cuspidal 
spectrum. 

In the meantime Sarnak's conjecture has been verified in  quite a number of 
cases. A. Reznikov proved it for congruence
groups in a group $G$ of real rank one, S. Miller \cite{Mi} proved  it
for $\G=\SL(3)$ and $\Gamma=\SL(3,\Z)$, the author \cite{Mu3} established it 
for $\G=\SL(n)$ and a congruence group $\Gamma$. The most general result
is due to Lindenstrauss and Venkatesh \cite{LV} who proved the following 
theorem.
\begin{theo}\label{th-weyllaw1}
Let $\G$ be a split adjoint semi-simple group over $\Q$ and let $\Gamma\subset
\G(\Q)$ be a congruence subgroup. Let $d=\dim S$. Then 
\begin{equation}\label{asymp3}
N_\Gamma^{\cu}(\lambda)\sim \frac{\vol(\Gamma\bs \widetilde X)}{(4\pi)^{d/2}
{\bf \Gamma}\left(\frac{d}{2}+1\right)}\lambda^{d/2},\quad \lambda\to\infty.
\end{equation}
\end{theo}
The method used by Lindenstrauss and Venkatesh is based on the construction of 
convolution operators with pure
cuspidal image. It avoids the delicate estimates of the contributions of the
Eisenstein series to the trace formula.
This proves existence of many cusp forms for these groups.

For an arbitrary $K$-type, we have the following theorem proved in \cite{Mu3}.
\begin{theo}\label{th-weyllaw2}
Let $n\ge 2$ and  $\widetilde X=\SL(n,\R)/\SO(n)$. Let 
$d=\dim\widetilde X=n(n+1)/2-1$. For
every principal congruence subgroup $\Gamma$ of $\SL(n,\Z)$ and every 
irreducible unitary representation $\sigma$ of $\SO(n)$ such that 
$\sigma|_{Z_\gamma}=\Id$, we have 
\begin{equation}\label{asymp4}
N_{\Gamma}^{\cu}(\lambda,\sigma)\sim\frac{\dim(\sigma)
\vol(\Gamma\bs \widetilde X)}{(4\pi)^{d/2}{\bf \Gamma}(d/2+1)}\lambda^{d/2}
\end{equation}
as $\lambda\to\infty$.
\end{theo}

The residual spectrum for $\SL(n)$ has been 
described by Moeglin and Waldspurger \cite{MW}. Combined with \eqref{limsup} it
follows that for $\G=\SL(n)$ we have
\begin{equation}\label{residual}
N_{\Gamma(N)}^{\res}(\lambda,\sigma)\ll \lambda^{d/2-1},
\end{equation}
where $d=\dim \SL(n,\R)/\SO(n)$ and $\Gamma(N)\subset \SL(n,\Z)$ is the 
principal congruence subgroup of level $N$. 

The proof of Theorem \ref{th-weyllaw2} uses the Arthur trace formula
combined with the heat equation method similar to the proof 
of \eqref{4.8}. The application of the Arthur trace formula requires the
adelic reformulation of the problem. 

We briefly describe the method. For all details we refer to \cite{Mu5}.
For simplicity we consider only the trivial $K_\infty$-type, i.e, we consider the
counting function $N_\Gamma^{\cu}(\lambda)$. 
By \eqref{residual} we can replace the counting function 
$N_\Gamma^{\cu}(\lambda)$ by $N_\Gamma^{\di}(\lambda)$. 
Let $\G=\GL(n)$ regarded as an algebraic group over $\Q$. Denote by $A_G$ 
the split component of the center of $\G$ and let $A_G(\R)^0$ be the 
component of 1 in $A_G(\R)$. Let $\Pi_{\di}(\G(\A),\xi_0)$ be the set of all 
irreducible subrepresentations of the regular representation of $G(\A)$ in
$L^2(\G(\Q)A_G(\R)^0\bs\G(\A))$. Given a representation 
$\pi\in\Pi_{\di}(\G(\A),\xi_0)$, let
$m(\pi)$ denote the multiplicity with which $\pi$ occurs in 
$L^2(\G(\Q)A_G(\R)^0\bs\G(\A))$. For any irreducible representation
$\pi=\pi_\infty\otimes\pi_f$ of $\G(\A)$, let $\H_{\pi_\infty}$ and $\H_{\pi_f}$
denote the Hilbert space of the representation $\pi_\infty$ and $\pi_f$,
respectively. Let $K_f$ be an open compact subgroup of $\G(\A_f)$. Denote
by $\H_{\pi_f}^{K_f}$ the subspace of $K_f$-invariant vectors in $\H_{\pi_f}$ and 
by $\H_{\pi_\infty}^{K_\infty}$ the subspace of $K_\infty$-invariant vectors in
$\H_{\pi_\infty}$. Given $\pi\in\Pi(\G(\A),\xi_0)$, denote by
$\lambda_{\pi_\infty}$ the Casimir eigenvalue of the restriction of $\pi_\infty$
to $\G(\R)^1$. Assume that $-1\neq K_f$. Then \eqref{asymp4} for the trivial
$K_\infty$-type follows by Karamata's theorem \cite[p. 446]{Fe} from the
 existence of an asymptotic expansion of the form
\begin{equation}\label{weyllaw4}
\sum_{\pi\in\Pi_{\di}(\G(\A),\xi_0)}
m(\pi)e^{t\lambda_{\pi_\infty}}\dim\bigl(\H_{\pi_f}^{K_f}\bigr)
\dim\bigl(\H_{\pi_\infty}^{K_\infty})\sim
\frac{\vol(\G(\Q)\bs \G(\A)^1/K_f)}{(4\pi)^{d/2}}
t^{-d/2}
\end{equation}
as $t\to +0$.

To establish \eqref{weyllaw4} we apply the Arthur trace formula as
follows. We choose a certain family of test functions $\widetilde\phi_t^1
\in C^\infty_c(\G(\A)^1)$, depending on $t>0$, which at the infinite place 
are given by the heat kernel $h_t\in C^\infty(\G(\R)^1)$ of the Laplacian 
$\widetilde\Delta$ on $\widetilde X$, multiplied by a certain cutoff function 
$\varphi_t$, and which
at the finite places is given by the normalized characteristic function of
an open compact subgroup $K_f$ of $\G(\A_f)$. Then by the non-invariant trace 
formula \cite{Ar1} we have the equality
\[
J_{\spec}(\widetilde\phi_t^1)=J_{\geo}(\widetilde\phi_t^1),\quad t>0.
\] 
Then we study asymptotic behavior of the spectral and the geometric side as 
$t\to0$. To deal with the geometric side, we use the fine
 ${\ho}$-expansion \cite{Ar6} 
\begin{equation}\label{fineexp}
J_{\geo}(f)=\sum_{\M\in\cL}\sum_{\gamma\in(\M(\Q_S))_{M,S}}a^M(S,\gamma)
J_M(\gamma,f),
\end{equation}
which expresses the distribution $J_{\geo}(f)$ in terms of weighted  orbital 
integrals $J_M(\gamma,f)$.
Here $\M$ runs over the set of Levi subgroups $\cL$ containing the Levi
component $\M_0$ of the standard minimal parabolic subgroup $\bP_0$, $S$ is
a finite set of places of $\Q$, and $(\M(\Q_S))_{M,S}$ is a certain set of
equivalence classes in $\M(\Q_S)$. This reduces our problem to the 
investigation of weighted orbital integrals. The key result is that 
$$\lim_{t\to0}t^{d/2}J_M(\widetilde\phi_t^1,\gamma)=0,$$
unless  $\M=\G$ and $\gamma =1$. This follows from the description of the local
weighted orbital integrals by \cite[Corollary 6.2]{Ar4}.
The contributions to \eqref{fineexp} of the terms where $\M=\G$ and
 $\gamma=1$ are easy to determine. Using
the behavior of the heat kernel $h_t(1)$ as $t\to0$, it follows that
\begin{equation}\label{asymp7}
J_{\geo}(\widetilde\phi_t^1)\sim\frac{\vol(\G(\Q)\bs \G(\A)^1/K_f)}{(4\pi)^{d/2}}
t^{-d/2}
\end{equation}
as $t\to0$. To deal with the spectral side we use Theorem \ref{thm-specexpand}.
This theorem allows us to replace $\widetilde\phi_t^1$ by a similar function
$\phi_t^1\in{\mathcal C}^1(G(\A)^1)$ which is given as the product of the 
heat kernel $h_t$ at infinity and the normalized characteristic function of
$K_f$. The term in $J_{\spec}(\phi^1_t)$ corresponding to $\M=\G$ is 
$J_{\spec,G}(\phi_t^1)=\tr R_{\di}(\phi_t^1)$, which is equal to the left hand 
side of \eqref{weyllaw4}. If $\M$ is a proper Levi subgroup of $\G$, then
$J_{\spec,M}(\phi_t^1)$ is given by \eqref{specside2}, which is a finite
some of integrals. The main ingredient of the integrals are logarithmic
derivatives of intertwining operators and the estimation of these integrals is
reduced to the estimation of the logarithmic derivatives.
Using  \eqref{normalization} this problem is reduced to the estimation of the
logarithmic derivatives of the normalizing factors and the local intertwining
operators. In the case of $\G=\GL(n)$, the normalizing factors are expressed in
terms of Ranking-Selberg $L$-functions \eqref{rankin-selb}. Using the 
analytic properties of Rankin-Selberg $L$-functions, it follows that there 
exist $C>0$ and $T>1$ such that for $\pi=\pi_1\otimes\pi_2$, 
$\pi_i\in\Pi_{\di}(\GL(n_i,\A))$,  we have
\begin{equation}\label{logderiv1}
\int_{T}^{T+1}\left|\frac{n_\alpha^\prime(\pi,i\lambda)}{n_\alpha(\pi,i\lambda)}
\right|\,d\lambda\le C\log(T+\nu(\pi_1\times \tilde\pi_2)),
\end{equation} 
where $\nu(\pi_1\times \tilde\pi_2)=N(\pi_1\times\tilde\pi_2)(2+c(\pi_1\times
\tilde\pi_2)$, $N(\pi_1\times\tilde\pi_2)$ is the conductor occurring in the 
functional equation \eqref{functequ} and  $c(\pi_1\times\tilde\pi_2)$ is the
analytic conductor defined in \cite[(4.21)]{Mu5}. For the proof of 
\eqref{logderiv1} see \cite[Proposition 5.1]{Mu5}. In the case of $\SL(2,\R)$
we have the pointwise estimate \eqref{4.6}. If we integrate it, we get the
analogue of \eqref{logderiv1} which would suffice to derive the Weyl law for
the principal congruence subgroups of $\SL(2,\Z)$. 
 
Finally we have to deal with normalized intertwining operators 
\[
R_{Q|P}(\pi,s)=\otimes_vR_{Q|P}(\pi_v,s).
\]
Since the open compact subgroup 
$K_{\fin}$ of $\G(\A_{\fin})$ is fixed, there are only finitely many places $v$
for which we have to  consider $R_{Q|P}(\pi_v,s)$. The main ingredient 
for the estimation of the logarithmic derivative of $R_{Q|P}(\pi_v,s)$, which
is uniform in $\pi_v$,  is a weak version of the Ramanujan conjecture 
(see \cite[Proposition 0.2]{MS}). 

Combining these estimations, it follows that for every proper Levi 
subgroup $\M$ of $\G$ we have
\begin{equation}\label{asymp9}
J_{\spec,M}(\phi_t^1)=O(t^{-(d-1)/2})
\end{equation}
as $t\to +0$. This proves \eqref{weyllaw4}. 

The next problem is to estimate the remainder term in the Weyl law. For 
$\G=\SL(n)$ this problem has been studied by E. Lapid and the author in 
\cite{LM}. Actually, we consider not only the cuspidal spectrum of the 
Laplacian, but the cuspidal spectrum  of the whole algebra of invariant 
differential operators $\cD(\widetilde X)$. 

As $\cD(\widetilde X)$ preserves the space of cusp
forms, we can proceed as in the compact case and decompose 
$L^2_{\cu}(\Gamma\bs \widetilde X)$ into joint eigenspaces of 
$\cD(\widetilde X)$. Recall that the characters of $\cD(\widetilde X)$ are
parametrized by $\af_\C^\ast/W$. Given $\lambda\in\af^\ast_\C/W$, denote by
$\chi_\lambda$ the corresponding character of $\cD(\widetilde X)$ and let
\[
\E_{\cu}(\lambda)=\left\{\varphi\in L^2_{\cu}(\Gamma\bs \widetilde X)\colon
D\varphi=\chi_\lambda(D)\varphi\right\}
\]
be the associated joint eigenspace.
Each eigenspace is finite-dimensional. Let $m(\lambda)=\dim\E_{\cu}(\lambda)$. 
Define the cuspidal spectrum 
$\Lambda_{\cu}(\Gamma)$ to be
\[\Lambda_{\cu}(\Gamma)=\{\lambda\in\af^\ast_\C/W\colon m(\lambda)>0\}.\]
Then we have an orthogonal direct sum decomposition
\[L^2_{\cu}(\Gamma\bs \widetilde X)=\bigoplus_{\lambda\in\Lambda_{\cu}(\Gamma)}
\E_{\cu}(\lambda).\]
Let $\beta(\lambda)$ be the Plancherel measure on $i\af^\ast$. 
Then in \cite{LM}  we established the following extension of 
main results of \cite{DKV} to congruence quotients of $S=\SL(n,\R)/\SO(n)$.
\begin{theo}\label{th-remainder}
Let  $d=\dim \widetilde X$. 
Let $\Omega\subset \af^\ast$ be a bounded domain with piecewise smooth
boundary. Then for $N\ge 3$ we have
\begin{equation}\label{5.5} 
\sum_{\substack{\lambda\in\Lambda_{\cu}(\Gamma(N))\\
\lambda\in i t\Omega}}m(\lambda)=
\frac{\vol(\Gamma(N)\bs \widetilde X)}{|W|}\int_{it\Omega}\beta(\lambda)\ 
d\lambda+
O\left(t^{d-1}(\log t)^{\max(n,3)}\right), 
\end{equation}
as $t\to\infty$, and
\begin{equation}\label{5.6} 
\sum_{\substack{\lambda\in\Lambda_{\cu}(\Gamma(N))\\
\lambda\in B_t(0)\setminus i\af^\ast}}m(\lambda)=
O\left(t^{d-2}\right),\quad t\to\infty.
\end{equation}
\end{theo}  
If we apply (\ref{5.5}) and (\ref{5.6}) to the unit ball in $\af^\ast$, we
get the following corollary.

\begin{corollary}\label{weyllaw5} 
Let $\widetilde X=\SL(n,\R)/\SO(n)$ and $d=\dim\widetilde X$.
Let $\Gamma(N)$ be the principal congruence subgroup of $\SL(n,\Z)$ of level 
$N$. Then for $N\ge 3$ we have
\[
N^{\cu}_{\Gamma(N)}(\lambda)=\frac{\vol(\Gamma(N)\bs \widetilde X)}
{(4\pi)^{d/2}{\bf\Gamma}\left(\frac{d}{2}+1\right)}\lambda^{d/2}+
O\left(\lambda^{(d-1)/2}(\log \lambda)^{\max(n,3)}\right),\quad \lambda\to\infty.
\]
\end{corollary} 
The condition $N\ge 3$ is imposed for technical reasons. It guarantees that
the principal congruence subgroup $\Gamma(N)$ is neat in the sense of Borel,
and in particular, has no torsion. This simplifies the analysis by eliminating
the contributions of the non-unipotent conjugacy classes in the trace formula.

Note that $\Lambda_{\cu}(\Gamma(N))\cap i\af^\ast$ is the cuspidal tempered
spherical spectrum. The Ramanujan conjecture \cite{Sa3}
 for $\GL(n)$ at the Archimedean place states that
\[\Lambda_{\cu}(\Gamma(N))\subset i\af^\ast\]
so that (\ref{5.6}) is empty, if the Ramanujan conjecture is true. However, 
the Ramanujan conjecture is far from 
being proved. Moreover, it is known to be false for other groups $\G$ and
(\ref{5.6}) is what one can expect in general.

The method to prove Theorem \ref{th-remainder} is an extension of the method of
\cite{DKV}. The Selberg trace formula, which is one of the basic tools in
 \cite{DKV}, is replaced by the non-invariant Arthur trace formula.  
Again, one of the main issues in the proof is the estimation of the logarithmic
derivatives of the intertwining operators occurring on the spectral side of the 
trace formula.

\subsection{Upper and lower bounds}

In some cases it suffices to have upper or lower bounds for the counting 
function. For example, Donnelly's result \eqref{limsup} implies that there
exists a constant $C>0$ such that
\begin{equation}
N_{\Gamma}^{\cu}(\lambda;\sigma)\le C(1+\lambda^{d/2}),\quad \lambda\ge 0.
\end{equation}
For the full discrete spectrum we have the bound \eqref{upperbd1}. However,
the exponent is not the optimal one. For some applications it is necessary
to have such a bound which is uniform in $\Gamma$. For the cuspidal spectrum
this problem has been studied by Deitmar and Hoffmann \cite{DH}. To state
the result, we have to introduce some notation. Let 
$\Gamma_n(N)$ be the principal congruence subgroup of $\GL(n,\Z)$ of level 
$N$. Let $\G$ be a connected reductive linear algebraic group over $\Q$. Let
$\eta\colon \G\to \GL(n)$ be a faithful $\Q$-rational representation. A family
${\mathcal T}$ of subgroups of $\G(\Q)$ is called a {\it family of  bounded 
depth} in $\G(\Q)$ if there exists $D\in\N$ which satisfies the following
property:  For every $\Gamma\in
{\mathcal T}$ there exists $N\in\N$ such that $\Gamma_n(N)\cap \eta(\G(\Q))$
is a subgroup of $\eta(\Gamma)$ of index at most $D$. We note that every
$\Gamma\in{\mathcal T}$ is contained in $\Gamma_0:=\Gamma_n(1)\cap\G(\Q)$.
Then the result of
Deitmar and Hoffmann \cite[Corollary 18]{DH} is the following theorem.
\begin{theo}
Let ${\mathcal T}$ be a family of bounded depth in $\G(\Q)$. There exists
$C>0$ such that for all $\Gamma\in{\mathcal T}$ and all $\lambda\ge 0$ we have
\begin{equation}\label{uniformbd}
N_\Gamma^{\cu}(\lambda;\sigma)\le C[\Gamma_0\colon\Gamma](1+\lambda)^{d/2}.
\end{equation}
\end{theo}
\begin{conjecture}\label{conj1}
The estimation \eqref{uniformbd} holds for $N_\Gamma^{\di}(\lambda;\sigma)$.
\end{conjecture}
Given the description of the residual spectrum for $\GL(n)$ by \cite{MW}, it
seems possible to establish this conjecture for $\GL(n)$. 

As for lower bounds there is the weak Weyl law established in \cite{LM}.
For $\sigma\in \widehat K$ let 
\[
c_\sigma(\Gamma)=\frac{\dim(\sigma)\vol(\Gamma\bs\widetilde X)}
{(4\pi)^{d/2}{\bf \Gamma}(d/2+1)}
\]
be the constant in Weyl's law, where $d=\dim(\widetilde X)$. 
Let $\G$ be a semisimple algebraic group defined over 
$\Q$ and let $\Gamma\subset \G(\Q)$ be a congruence subgroup defined by an
open compact subgroup $K_{\fin}=\prod_pK_p$ of $\G(\A_{\fin})$.   Let $S$ be a
finite set of primes. We will say that $\Gamma$ is deep enough with respect
to $S$, if for every prime $p\in S$, $K_p$ is a subgroup of some minimal 
parahoric subgroup of $\G(\Q_p)$. Then the main result of \cite{LM} is the
following theorem.
\begin{theo}
Let $\G$ be an almost simple connected and simply connected semisimple
algebraic group defined over $\Q$ such that $\G(\R)$ is non compact. 
Let $S$ be a finite set of primes containing at least two primes. Then for
every congruence subgroup $\Gamma\subset \G(\Q)$ there exists a nonnegative
constant $c_S(\Gamma)\le 1$ such that for every $\sigma\in\widehat K$ with
$\sigma|_{Z_\Gamma}=\Id$ we have
\[
c_\sigma(\Gamma)c_S(\Gamma)\le \liminf_{\lambda\to\infty}
\frac{N_\Gamma^{\cu}(\lambda,\sigma)}{\lambda^{d/2}}.
\]
Moreover $c_S(\Gamma)>0$ if $\Gamma$ is deep enough with respect to $S$.
\end{theo}

\subsection{Self-dual automorphic representations}

So far, we considered only the family of all cusp forms of $\GL(n,\A)$. 
A nontrivial subfamily is formed by the family of self-dual automorphic 
representations. They arise as functorial lifts of automorphic representations
of classical groups. Functoriality from quasisplit classical groups to 
general linear groups has been
established by Cogdell, Kim, Piatetski-Shapiro, and Shahidi for generic
automorphic representations and then by Arthur for all representations.
In his thesis, V. Kala has studied the counting function of self-dual
cuspidal automorphic representations of $\GL(n,\A)$. For $N\in \N$ with
prime decomposition $N=\prod_p p^{r(p)}$ let
\[
K_p(N):=\left\{k\in\GL(n,\Z_p)\colon k\equiv 1\mod p^{r(p)}\Z_p\right\}
\]
Let $K(N)$ be the principal congruence subgroup defined by 
\[
K(N):=O(n)\times\prod_p K_p(N).
\]
Let
\[
N_{\sd}^{K(N)}(\lambda):=\sum_{\substack{\lambda(\Pi)\le\lambda\\
\Pi\cong\widetilde\Pi}}\dim\Pi^{K(N)},
\]
where the sum ranges over all self-dual cuspidal automorphic representations
$\Pi$ of $\GL(n,\A)$ with Casimir eigenvalues $\le \lambda$. Then the main 
result of \cite{Ka} is the following theorem.

\begin{theo}\label{th-sd}
Let $n=2m+\varepsilon$ with $\varepsilon=0,1$. Put $d=m^2+m$.
For all $N\in\N$ there exist constants $C_1,C_2>0$ such that for $\lambda\gg 0$ 
one has
\[
C_1\lambda^{d/2}\le N_{\sd}^{K(N)}(\lambda)\le C_2\lambda^{d/2}.
\]
\end{theo}

By Corollary \ref{weyllaw5}, the counting function of all cuspidal 
representations, counted similarly, is asymptotic to $C\lambda^{d/2}$, where
$d=(n^2+n-2)/2$. Hence for $n>2$, the density of self-dual cusp forms is
zero. 

The main idea of the proof of Theorem \ref{th-sd} is to consider the 
descent $\pi$ of each self-dual
cuspidal automorphic representation $\Pi$ of $\GL(n,\A)$ to one of the 
quasisplit classical groups $\G(\A)$ and to use results towards the Weyl
law on $\G(\A)$. The number $d=m^2+m$ is related to the dimension of the
corresponding symmetric space $\G(\R)/K_\infty$ (see \cite[p.17]{Ka}).  
The key problem of the proof is to relate the Casimir eigenvalue and the
existence of $K(N)$-fixed vectors for $\Pi$ and $\pi$. 

In a special case Kala's method leads to an exact asymptotic formula. Let
$n=2m$ and $d=m^2+m$. Let $K=O(n)\times \prod_p K_p$ with $K_p=\GL(n,\Z_p)$. 
Then there exists $C>0$ such that
\begin{equation}
N_{\sd}^K(\lambda)=C\lambda^{d/2}+o(\lambda^{d/2})
\end{equation}
(see \cite[Corollary 6.2.2]{Ka}). One may conjecture that this is true in 
general.

\subsection{Weyl's law for Hecke operators}

One can also study the asymptotic distribution of infinitesimal characters of
cuspidal automorphic representations weighted by the eigenvalues of Hecke
operators acting on cusp forms of $\GL(n)$. For details we refer to the recent
papers by J. Matz \cite{Ma1}, J. Matz and N. Templier \cite{MT} and the
survey article of J. Matz in these proceedings.

\section{The limit multiplicity problem} 
\setcounter{equation}{0}

The limit multiplicity problem is another basic problem which is concerned with
the asymptotic behavior of automorphic spectra.

In this section we summarize some of the known results about the limit
multiplicity problem. To begin with we recall
some facts concerning the Plancherel measure $\mu_{\pl}$ on $\Pi(G)$.
First of all, the support of $\mu_{\pl}$ is the tempered dual $\Pi(G)_{\temp}$,
consisting of the equivalence classes of the irreducible unitary tempered 
representations. Up to a closed subset of Plancherel measure zero, the 
topological space $\Pi(G)_{\temp}$ is homeomorphic to a countable union of
Euclidean spaces of bounded dimensions. Under this homeomorphism the Plancherel
density is given by a continuous function. We call the relatively quasi-compact
subsets of $\Pi(G)$ {\it bounded}. We note that $\mu_\Gamma(A)<\infty$ for
bounded sets $A\subset \Pi(G)$ under the reduction-theoretic assumptions on
$(G,\Gamma)$ mentioned above (see \cite{BG}). A bounded subset $A$ of 
$\Pi(G)_{\temp}$ is called a Jordan measurable subset, if $\mu_{\pl}(\partial A)
=0$, where $\partial A=\bar A- \inter(A)$ is the boundary of $A$ in 
$\Pi(G)_{\temp}$. Furthermore, a Riemann integrable function on $\Pi(G)_{\temp}$
is a bounded, compactly supported function which is continuous almost 
everywhere with respect to the Plancherel measure. 

Let $(\mu_n)_{n\in\N}$ be a sequence of Borel measures on $\Pi(G)$. We say that 
the sequence $(\mu_n)_{n\in\N}$ has the {\it limit multiplicity property}
(property (LM)), if the following two conditions are satisfied. 
\begin{enumerate}
\item[1)] For every Jordan measurable set $A\subset \Pi(G)_{\temp}$ we have
\[
\mu_n(A)\to \mu_{\pl}(A),\quad \text{as}\;\; n\to\infty.
\]
\item[2)] For every bounded subset $A\subset \Pi(G)\setminus \Pi(G)_{\temp}$
we have 
\[
\mu_n(A)\to 0, \quad \text{as}\;\; n\to\infty.
\]
\end{enumerate}
We note that condition 1) can be restated as 
\begin{enumerate}
\item[1a)] For every Riemann integrable function $f$ on $\Pi(G)_{\temp}$ one has
\[
\lim_{n\to\infty}\mu_{n}(f)=\mu_{\pl}(f).
\]
\end{enumerate}
Now let $(\Gamma_n)_{n\in\N}$ be a sequence of lattices in $G$. The sequence
$(\Gamma_n)_{n\in\N}$ is said to have the limit multiplicity property (LM), if
the sequence of measures $(\mu_{\Gamma_n})_{n\in\N}$ has property (LM).  

The limit multiplicity problem can be formulated as follows: under which 
conditions does the sequence of measures $\mu_{\Gamma_n}$ satisfy property (LM)?

The limit multiplicity problem has been studied to a great extent in the case of
uniform lattices. In this case, $R_\Gamma$ decomposes discretely. It started with
the work of DeGeorge and Wallach \cite{DW1},\cite{DW2}, who 
considered towers of normal subgroups, i.e.,  descending
sequences of normal subgroups of finite index of a given uniform lattice with
trivial intersection. 
For such sequences they  dealt with
the case of discrete series representations and the tempered spectrum, if the
split rank of $G$ is $1$. Subsequently, Delorme \cite{De} solved 
the limit multiplicity problem affirmatively for normal towers of cocompact 
lattices. 
Recently, there has been great progress in proving limit multiplicity for much 
more general
sequences of uniform lattices by Abert et al \cite{AB1},\cite{AB2}. In 
particular, families of 
non-commensurable lattices were considered for the first time. The basic idea
is the notion of Benjamini-Schramm convergence (BS-convergence), which 
originally was introduced for sequences of finite graphs of bounded degree
and has been adopted by Abert et al to sequences of Riemannian manifolds. For
a Riemannian manifold $M$ and $R>0$ let 
\[
M_{<R}=\{x\in M\colon \injrad_M(x)<R\}.
\]
Let $(\Gamma_n)$ be a sequence of lattices in $G$. Then the orbifolds 
$M_n=\Gamma_n\bs X$ are said to {\it BS-converge} to $X$, if for every $R>0$
one has
\begin{equation}\label{bs-conv}
\lim_{n\to+\infty}\frac{\vol((M_n)_{<R})}{\vol(M_n)}=0.
\end{equation}
To find examples of sequences $(\Gamma_n)$ which satisfy this condition,
consider a cocompact arithmetic lattice $\Gamma_0\subset G$. By
\cite[Theorem 5.2]{AB1} there exist constants $c,\mu>0$ such that for any
congruence subgroup $\Gamma\subset \Gamma_0$ and any $R>1$ one has
\begin{equation}\label{bs-conv1}
\vol((\Gamma\bs X)_{<R})\le e^{cR} \vol(\Gamma\bs X)^{1-\mu}.
\end{equation}
Thus any sequence $(\Gamma_n)$ of congruences subgroups of $\Gamma_0$ such
that $\vol(\Gamma_n\bs G)\to \infty$ as $n\to\infty$ satisfies \eqref{bs-conv}.

A family of lattices in $G$ is called to be uniformly discrete, if there 
exists a neighborhood of the identity in $G$ that intersects trivially all of
their conjugates. For torsion-free lattices $\Gamma_n$ this is equivalent to
the condition that there is a uniform lower bound of the injectivity radii of
the manifolds $\Gamma_n\bs X$. In particular, any family of normal subgroups
 $(\Gamma_n)$ of a fixed uniform lattice $\Gamma$ is uniformly discrete.
Now the following theorem is one of the main results of 
\cite[Theorem 1.2]{AB1}.
\begin{theo}[\cite{AB1}]
Let $(\Gamma_n)$ be a uniformly discrete sequence of lattices in $G$ such that
the orbifolds $\Gamma_n\bs X$ BS-converge to $X$. Then the sequence 
$(\Gamma_n)$ has the (LM) property.
\end{theo}
It follows from the discussions above that any sequence of congruence subgroups
$(\Gamma_n)$ of a given cocompact arithmetic lattice $\Gamma_0$ of $G$
satisfies the assumptions of the theorem. 

A special case of the limit multiplicity property is the case of a singleton
$A=\{\pi\}$. Let $\Pi(G)_d\subset \Pi(G)$ be the discrete series and 
$d(\pi)$ the formal degree of $\pi\in\Pi(G)_d$. If $(\Gamma_n)$ is a sequence
of lattices in $G$ which satisfies the (LM) property, then it follows that
\begin{equation}\label{discr-ser}
\lim_{n\to\infty}\frac{m_{\Gamma_n}(\pi)}{\vol(\Gamma_n\bs G)}=\begin{cases} d(\pi),
&\pi\in\Pi(G)_d,\\ 0, &\text{else}.
\end{cases}
\end{equation}
It was first proved by DeGeorge and Wallach \cite{DW1} that \eqref{discr-ser}
holds for any tower of normal subgroups of a given uniform lattice of $G$. 

An important problem is to extend these results to the non-cocompact case.
Then the spectrum contains a continuous part and much less is known. The limit
multiplicity problem has been solved for normal towers of arithmetic lattices
and discrete series $L$-packets of representations (with regular parameters)
by Rohlfs and Speh \cite{RoS}. Then Savin \cite{Sav} solved the limit
multiplicity problem for the discrete series and
 normal towers of congruence subgroups. 

In \cite{FLM2} we dealt with the general case. Let $F$ be a number field and
denote by $\cO_F$ its ring of integers. For the non-compact lattice
$\SL(n,\cO_F)\subset \SL(n,F\otimes \R)$ we have the following 
result.
\begin{theo}\label{lm1}
Let $F$ be a number field. 
Then the collection of principal congruence subgroups $(\Gamma_N)$ of 
$\SL(n,{\mathcal O}_F)$ has the limit multiplicity property.
\end{theo}
In \cite{FL2}, T. Finis and E. Lapid extended this result to the collection of 
all
congruence subgroups of $\SL(n,{\mathcal O}_F)$, not containing non-trivial
central elements. In \cite{FLM2}, we also discussed the case of a general 
reductive group. 

\subsection{The density principle and the trace formula}\label{density} 

A standard approach to the limit multiplicity problem is to use integration
against test functions on $G$ and the trace formula. Let $K$ be a maximal
compact subgroup of $G$. Denote by $C^\infty_{c,\fin}(G)$ the space of smooth,
compactly supported bi-$K$-finite functions on $G$. Given $f\in 
C^\infty_{c,\fin}(G)$, define $\hat f(\pi)$ for $\pi\in\Pi(G)$ by $\hat f(\pi):=
\tr \pi(f)$. The function $\pi\in\Pi(G)\mapsto \hat f(\pi)$ on
$\Pi(G)$  is the ``Fourier transform'' of $f$.  Let $\mu$ be a Borel measure on
$\Pi(G)$. Then $\mu(\hat f)$ is defined (of course, it might be divergent). 
In particular, we have the two Borel measures $\mu_{\pl}$ and $\mu_\Gamma$ 
defined on $\Pi(G)$. For these measures we have $\mu_{pl}(\hat f)=f(1)$ and
\begin{equation}
\mu_\Gamma(\hat f)=\frac{1}{\vol(\Gamma\bs G)}\tr R_{\Gamma,\di}(f).
\end{equation}
By \cite{Mu2}, $R_{\Gamma,\di}(f)$ is a trace class operator. Thus the right hand 
side is well defined. Furthermore, by the Plancherel theorem we have
$\mu_{\pl}(\hat f)=f(1)$. The density principle of Sauvageot \cite{Sau}, which
is a refinement of the work of Delorme, can be stated as follows.

\begin{theo}\label{density-princ}
Let $(\mu_n)_{n\in\N}$ be a sequence of Borel measures on $\Pi(G)$ and assume that
for all $f\in C^\infty_{c,\fin}(G)$ we have 
\begin{equation}\label{conver1}
\mu_n(\hat f)\to \mu_{\pl}(\hat f)=f(1), \quad \text{as}\;\;n\to\infty.
\end{equation}
Then $(\mu_n)_{n\in\N}$ satisfies (LM).
\end{theo}

Now let $(\Gamma_n)_{n\in\N}$ be a sequence of lattices in $G$. Then by Theorem
\ref{density-princ} it follows that $(\Gamma_n)_{n\in\N}$ satisfies (LM), if
\begin{equation}\label{conv-measure}
\mu_{\Gamma_n}(\hat f)\to f(1),\quad n\to\infty,
\end{equation}
for all $f\in C^\infty_{c,\fin}(G)$.  
A standard approach to verify \eqref{conv-measure} is to use the trace
formula. In the case of co-compact lattices this is rather simple. Let $\Gamma$
be a cocompact lattice in $G$. Then the Selberg trace formula is the following
equality
\[
\vol(\Gamma\bs G)\mu_\Gamma(\hat f)=\tr R_\Gamma(f)=\sum_{\{\gamma\}\in C(\Gamma)}
\vol(\Gamma_\gamma\bs G_\gamma)\int_{G_\gamma\bs G} f(x^{-1}\gamma x)\;dx,
\]
where  $C(\Gamma)$ denotes the $\Gamma$-conjugacy classes  of 
$\Gamma$,
and $G_\gamma$ (resp. $\Gamma_\gamma$) denotes the centralizer of $\gamma$ in
$G$ (resp. $\Gamma$). Let $\Gamma_1\subset\Gamma$ be a finite index subgroup.
For $\gamma\in\Gamma$ let
\begin{equation}\label{const1}
c_{\Gamma_1}(\gamma)=|\{\delta\in\Gamma_1\bs\Gamma\colon \delta\gamma\delta^{-1}
\in\Gamma_1\}|.
\end{equation}
In \cite{Co}, Corwin shows that the elements on the right hand side of the
trace formula for $\Gamma_1$ can be grouped together in a way to give 
\begin{equation}\label{corwin1}
\mu_{\Gamma_1}(\hat f)=\frac{1}{\vol(\Gamma\bs G)}\sum_{\{\gamma\}\in C(\Gamma)}
\vol(\Gamma_\gamma\bs G_\gamma)
\frac{c_{\Gamma_1}(\gamma)}{[\Gamma\colon\Gamma_1]}\int_{G_\gamma\bs G}
f(x^{-1}\gamma x)\;dx.
\end{equation}
For a central element $\gamma$ we obviously have $c_{\Gamma_1}(\gamma)
=[\Gamma\colon\Gamma_1]$. Assume that the center of $\Gamma$ is trivial. Let
$(\Gamma_n)_{n\in\N}$ be a sequence of finite index subgroups of $\Gamma$. Then
we have
\begin{equation}\label{corwin2}
\mu_{\Gamma_n}(\hat f)=f(1)+\frac{1}{\vol(\Gamma\bs G)}
\sum_{\{\gamma\}\in C(\Gamma)\setminus \{1\}}\vol(\Gamma_\gamma\bs G_\gamma)
\frac{c_{\Gamma_n}(\gamma)}{[\Gamma\colon\Gamma_n]}\int_{G_\gamma\bs G}
f(x^{-1}\gamma x)\;dx.
\end{equation}
By dominated convergence, it follows that in order to establish 
\eqref{conver1} for the sequence $(\Gamma_n)_{n\in\N}$,
it suffices to show that for every $\gamma\in\Gamma$, $\gamma\neq 1$, we have
\begin{equation}\label{converg2}
\frac{c_{\Gamma_n}(\gamma)}{[\Gamma\colon \Gamma_n]}\to 0,\quad\text{as}\;\;
n\to\infty.
\end{equation}
Now note that if $\Gamma_1$ is a normal subgroup of $\Gamma$, then 
$c_{\Gamma_1}(\gamma)/[\Gamma\colon\Gamma_1]$ is the characteristic function of 
$\Gamma_1$. Thus for normal towers of finite index subgroups of $\Gamma$ the
condition \eqref{converg2} holds trivially. This implies Delorme's result.
 
If $\Gamma$ is not co-compact, the Selberg trace formula is only available in
the rank one case. We have to switch to the adelic framework so that we can use
the Arthur trace formula. 

Thus let now $\G$ be an arbitrary reductive group defined over $\Q$. Let
$\A =\R\times \A_{\fin}$ be the locally compact adele ring of $\Q$. For every 
place $v$ of $\Q$ (i.e. $v=\infty$ or $v=p$ a prime) let $|\cdot|_v$ be the
normalized absolute value of $\Q$. As usual. $\G(\R)^1$ denotes the 
intersection of the kernels of the homomorphisms $|\chi|\colon \G(\R)\to \R^+$,
where $\chi$ runs over the $\Q$-rational characters of $\G$. Similarly we define
the normal subgroup $\G(\A)^1$ of $\G(\A)$. Every $\pi\in\Pi(\G(\A)^1)$ can be 
written as $\pi=\pi_\infty\otimes \pi_{\fin}$, where $\pi_\infty\in\Pi(\G(\R)^1)$
and $\pi_{\fin}\in\Pi(\G(\A_{\fin}))$. 
Fix a Haar measure on $\G(\A)$. For
any open compact subgroup $K_f$ of $\G(\A_{\fin})$, let $\mu_K=\mu_K^G$ be the
measure on $\Pi(\G(\R)^1)$ defined by
\begin{equation}\label{adelic-meas}
\begin{split}
\mu_K&=\frac{1}{\vol(\G(\Q)\bs \G(\A)^1/K)}\sum_{\pi\in\Pi(G(\R)^1)}
\Hom_{\G(\R)^1}(\pi,L^2(\G(\Q)\bs \G(\A)^1/K))\delta_\pi\\
&=\frac{\vol(K)}{\vol(\G(\Q)\bs \G(\A)^1)}\sum_{\pi\in\Pi(\G(\A)^1)}
\dim\Hom_{\G(\A)^1}(\pi,L^2(\G(\Q)\bs \G(\A)^1))\dim(\pi_{\fin})^K\delta_{\pi_\infty}.
\end{split}
\end{equation}
We say that a sequence $(K_n)_{n\in\N}$ of open compact subgroups of 
$\G(\A_{\fin})$
has the limit multiplicity property, if $\mu_{K_n}\to\mu_{\pl}$, $n\to\infty$, in
the sense that
\begin{enumerate}
\item[(1)] For every Jordan measurable subset $A\subset \Pi(\G(\R)^1))_{\temp}$
we have $\mu_{K_n}(A)\to \mu_{\pl}(A)$ as $n\to\infty$, and
\item[(2)] For every bounded subset $A\subset\Pi(\G(\R)^1)\setminus 
\Pi(\G(\R)^1))_{\temp}$, we have $\mu_{K_n}(A)\to 0$ as $n\to\infty$.
\end{enumerate}
Again we can rephrase the first condition by saying that for any Riemann
integrable function $f$ on $\Pi(\G(\R)^1)_{\temp}$ we have
\begin{equation}\label{riem-integr}
\mu_{K_n}(f)\to\mu_{\pl}(f),\quad\text{as}\;\;n\to\infty.
\end{equation}
Note that when $\G$ satisfies the strong approximation property (which is the 
case if $\G$ is semisimple, simply connected, and without any $\Q$-simple factor
${\bf H}$ for which ${\bf H}(\R)$ is compact) and $K$ is an open compact 
subgroup of $\G(\A_{\fin})$, then we have
\[
\G(\Q)\bs \G(\A)/K\cong \Gamma_K\bs \G(\R),
\]
where $\Gamma_K=\G(\Q)\cap K$ is a lattice in the connected semisimple Lie 
group $\G(\R)$. 

Now for $f\in C^\infty_{c,\fin}(\G(\R)^1)$ we have
\begin{equation}\label{measure4}
\mu_K(\hat f)=\frac{1}{\vol(\G(\Q)\bs \G(A)^1)}\tr R_{\di}(f\otimes {\bf 1}_K)
\end{equation}
and
\begin{equation}\label{measure5}
\mu_{\pl}(\hat f)=f(1).
\end{equation}
Sauvageot's density principle \cite{Sau} can now be reformulated as follows.
\begin{theo}
Let $(K_n)_{n\in\N}$ be a sequence of open compact subgroups of $\G(\A_{\fin})$. 
Suppose that for every $f\in C^\infty_{c,\fin}(\G(\R)^1)$ we have
\begin{equation}\label{converg3}
\mu_{K_n}(\hat f)\to f(1),\quad n\to \infty.
\end{equation}
Then $(K_n)_{n\in\N}$ has the limit multiplicity property.
\end{theo}
To try to verify \eqref{converg3}, it is natural to use Arthur's 
(non-invariant) trace formula, which is an equality 
\[
J_{\spec}(h)=J_{\geo}(h),\quad h\in C_c^\infty(\G(\A)^1),
\]
of two distribution on $\G(\A)^1$ \cite{Ar1}, \cite{Ar2}, \cite{Ar3}. 
The distribution $J_{\spec}$ is expressed in terms of spectral data and 
$J_{\geo}$ in terms of geometric data.  The main terms on the geometric side are
the elliptic orbital integrals. In particular, the contribution 
$\vol(\G(\Q)\bs \G(\A)^1)h(1)$ of the identity element occurs on the 
geometric side. The main term on the spectral side is $\tr R_{\di}(h)$. 
By \eqref{measure4} it follows that
\eqref{converg3} can be broken down into the following two statements. For
every $f\in C^\infty_{c,\fin}(G(\R)^1)$ we have
\begin{equation}\label{spec-conv}
J_{\spec}(f\otimes {\bf 1}_{K_n})-\tr R_{\di}(f\otimes{\bf 1}_{K_n})
\to 0, \quad n\to\infty,
\end{equation}
and
\begin{equation}\label{geo-conv}
J_{\geo}(f\otimes {\bf 1}_{K_n})\to \vol(\G(\Q)\bs\G(\A)^1)f(1),\quad
n\to\infty.
\end{equation}
We call \eqref{spec-conv} the spectral - and \eqref{geo-conv} the geometric
limit property.

\subsection{Bounds on co-rank one intertwining operators} 
\label{SectionBounds}

In this section we formulate two conditions on the behavior of
the intertwining operators $M_{Q|P}$ which imply the spectral limit property 
for a given $\G$. They also imply 
Weyl's law for the group $\G$. We call these properties \TWN\ (tempered 
winding number) and \BD\ (bounded degree). The first property is global and
second local. The first property is connected with analytic problems in the
theory of automorphic $L$-functions. 

We will use the notation $A\ll B$ to mean that there exists a constant $c$ 
(independent of the parameters under consideration) such that $A\le cB$.
If $c$ depends on some parameters (say $F$) and not on others then we will 
write $A\ll_FB$.

Fix a faithful $\Q$-rational representation $\rho: \G \to \GL (V)$ and a 
$\Z$-lattice $\Lambda$ in the representation space
$V$ such that the stabilizer of $\hat{\Lambda} = \hat\Z \otimes 
\Lambda \subset \A_{\fin} \otimes V$ in $\G(\A_{\fin})$ is the group $\K_{\fin}$.
(Since the maximal compact subgroups of $\GL (\A_{\fin} \otimes V)$ are 
precisely the stabilizers of lattices, it is easy to see that such a lattice 
exists.)
For any $N\in\N$ let
\begin{equation}\label{congruence}
\K (N) = \{ g \in \G (\A_{\fin}) \, : \, \rho (g) v \equiv v 
\pmod{N \hat{\Lambda}}, \quad v \in \hat{\Lambda} \}
\end{equation}
be the principal congruence subgroup of level $N$,
an open normal subgroup of $\K_{\fin}$. 
The groups $\K (N)$ form a neighborhood basis of the identity element in 
$\G(\A_{\fin})$. For an open subgroup $K$ of $\K_{\fin}$ let the \emph{level} of 
$K$ be the smallest integer $N$ such that $\K(N)\subset K$. 
Analogously, define $\level (K_v)$ for open subgroups $K_v \subset \K_v$.

As in \cite{Mu6}, for any $\pi\in\Pi(\M(\R))$ we define
$\param_\pi=\sqrt{\lambda_\pi^2+\lambda_\tau^2}$,
where $\tau$ is a lowest $\K_\infty$-type of $\Ind^{\G(\R)}_{\bP(\R)}(\pi)$
and $\lambda_\pi$ and $\lambda_\tau$ are the Casimir eigenvalues of $\pi$ and
$\tau$, respectively. Note that 
this is well-defined, because $\lambda_\tau$ is independent of $\tau$.
Roughly speaking, $\param_\pi$ measures the size of $\pi$.
For $\M \in\levis$, $\alpha \in \Sigma_M$ and $\pi\in\Pi_{\disc}(\M(\A))$ let
$n_\alpha (\pi,s)$ be the global normalizing factor defined by 
\eqref{normalization}.

\begin{definition}
We say that the group $\G$ satisfies the property \TWN\ 
(tempered winding number) if
for any $\M\in\levis$, $\M \neq \G$, and any finite subset $\types \subset 
\Pi (\K_{M,\infty})$ there exists an integer $k > 1$ such that for any 
$\alpha \in \Sigma_M$ and any $\epsilon>0$ we have
\begin{equation}\label{logderiv2}
\int_{\iii\R}\left| {\frac{n'_\alpha (\pi,s)}{n_\alpha(\pi,s)}} \right| 
(1+\abs{s})^{-k}\ ds\ll_{\types,\epsilon} (1+\param_{\pi_\infty})^k\level(K_M)^{\epsilon}
\end{equation}
for all open compact subgroups $K_M$ of $\K_{M,\fin}$ and all 
$\pi=\pi_\infty\otimes\pi_{\fin}\in\Pi_{\disc}(\M(\A))$ such that $\pi_\infty$ 
contains a $\K_{M,\infty}$-type in the set $\types$ and $\pi_{\fin}^{K_M}\neq 0$.
\end{definition}

Since the normalizing factors $n_\alpha (\pi,s)$ arise from co-rank one 
situations, the property \TWN\ is hereditary for Levi subgroups.

\begin{remark}
If we fix an open compact subgroup $K_M$, then the corresponding bound
\[
\int_{\iii\R}\left| {\frac{n'_\alpha (\pi,s)}{n_\alpha(\pi,s)}} \right| (1+\abs{s})^{-k}\ ds\ll_{K_M} (1+\param_{\pi_\infty})^k
\]
is the content of \cite[Theorem 5.3]{Mu6}. So, the
point of \TWN\ lies in the dependence of the bound on $K_M$.
\end{remark}

\begin{remark}
In fact, we expect that
\begin{equation}\label{logbnd}
\int_{T}^{T+1}\left| {\frac{n'_\alpha(\pi,\iii t)}{n_\alpha(\pi,\iii t)}} \right| \ dt\ll 1+\log(1+T)+\log(1+\param_{\pi_\infty})+\log\level(K_M)
\end{equation}
for all $T\in\R$ and $\pi\in\Pi_{\disc}(\M(\A))^{K_M}$.
This would give the following strengthening  of \TWN:
\[
\int_{\iii\R}\left| {\frac{n'_\alpha (\pi,s)}{n_\alpha(\pi,s)}} \right| (1+\abs{s})^{-2}\ ds\ll 1+\log(1+\param_{\pi_\infty})+\log\level(K_M)
\]
for any $\pi\in\Pi_{\disc}(\M(\A))^{K_M}$.
\end{remark}

\begin{remark} \label{GlobalRemark}
If $G'$ is simply connected, then by \cite[Lemma 1.6]{Lub} (cf.~also \cite[Proposition 1]{FLM2}) we can
replace $\level (K_M)$ by $\vol (K_M)^{-1}$ in the definition of \TWN\ (as well as in \eqref{logbnd}).
\end{remark}

For $\GL(n)$ the the normalizing factors are expressed in terms of 
Rankin-Selberg $L$ functions (see \eqref{rankin-selb}). The known properties 
of Rankin-Selberg $L$-functions lead to the estimation \eqref{logderiv1},
which implies the desired estimation. By \cite[Lemma 5.4]{FLM2}, the case
of $\SL(n)$ can be reduced to $\GL(n)$. In this way we get (see \cite{FLM2})

\begin{theo} \label{mainglobal}
The estimate \eqref{logbnd} holds for $\G=\GL(n)$ or $\SL(n)$ with an 
implied constant depending only on $n$.
In particular, the groups $\GL(n)$ and $\SL(n)$ satisfy the property \TWN.
\end{theo}

\begin{remark}
For general groups $\G$ the normalizing factors are given, at least up to 
local factors, by quotients of automorphic $L$-functions
associated to the irreducible constituents of the adjoint action of the 
$L$-group $^LM$ of $\M$ on the unipotent radical
of the corresponding parabolic subgroup of $^LG$ \cite{La}.
To argue as above, we would need to know that
these $L$-functions have finitely many poles and satisfy a functional 
equation with the associated conductor
bounded by an arbitrary power of $\level (K_M)$ for automorphic 
representations $\pi\in\Pi_{\disc}(\M(\A))^{K_M}$.
Unfortunately, finiteness of poles and the expected functional equation are 
not known in general.
It is possible that for classical groups these properties are within reach.
\end{remark} 

Now we come to the second condition, which is a condition on the local
intertwining operators. Recall that for a finite prime $p$, the matrix 
coefficients of the local 
normalized intertwining operators $R_{Q|P}(\pi_p,s)^{K_p}$ are rational functions
of $p^s$. Moreover, their denominators can be controlled in terms of $\pi_p$,
and the degrees of these denominators are bounded in terms of $\G$ only. For
any Levi subgroup $\M\in\cL$ let $\G_{\M}$ be the closed subgroup of $\G$ 
generated by the unipotent radicals $\bU_P$, where $\bP\in\cP(\M)$. It is a 
connected semisimple normal subgroup of $\G$. 

\begin{definition}
We say that $\G$ satisfies \BD\ (bounded degree) if there exists a constant $c$
(depending only on $\G$ and $\rho$), such that for any $\M\in\cL$, $\M\neq\G$,
and adjacent parabolic groups $\bP,\bQ\in\cP(\M)$, any prime $p$, any open 
subgroup $K_p\subset \K_p$ and any smooth irreducible representation $\pi_p$
of $\M(\Q_p)$, the degrees of the numerators of the linear operators
$R_{Q|P}(\pi_p,s)^{K_p}$ are bounded by $c\log_p\level^{G_M}(K_P)$ if $\K_p$ is 
hyperspecial, and by $c(1+\log_p\level^{\G_{\M}}(K_p))$, otherwise.
\end{definition}

Property \BD\ has been studied in \cite{FLM3}. By 
\cite[Theorem 1, Proposition 6]{FLM3} we have the following theorem.
\begin{theo}\label{thm-bd}
The groups $\GL(n)$ and $\SL(n)$ satisfy (BD).
\end{theo}
The property \BD\ has the following consequence.
\begin{prop}
Suppose that $\G$ satisfies (BD). Let $\M\in\cL$ and let $\bP,\bQ\in\cP(\M)$ 
be adjacent parabolic subgroups. Then for all $\pi\in\Pi_{\di}(\M(\A))$, for 
all open subgroups $K\subset \K_{\fin}$ and all $\tau\in\Pi(\K_\infty)$ we have
\begin{equation}\label{logderiv3}
\begin{split}
\int_{i\R}\biggl\|R_{Q|P}(\pi,s)^{-1}\frac{d}{ds}R_{Q|P}(\pi,s)
\Big|_{I_P^G(\pi)^{\tau,K}}&\biggr\|(1+|s|^2)^{-1}\;ds\\
&\ll 1+\log(\|\tau\|+\level(K;\G_{\M}^+)).
\end{split}
\end{equation}
\end{prop}
The proof of the proposition follows from a generalization of Bernstein's 
inequality \cite{BE}. Suppose that $\G$ satisfies (TWN) and (BD).
Combining \eqref{logderiv2} and \eqref{logderiv3} we get an appropriate 
estimate for the corresponding integral involving the logarithmic derivative
of the intertwining operators.

\subsection{Application to the limit multiplicity problem}

The limit multiplicity property is a consequence of properties (TWN) and
(BD). The proof proceeds by induction over the Levi subgroups of $\G$. 
The property that is suitable for the induction procedure is not the spectral
limit property, but a property that we call {\it polynomial boundedness}
(PB). This is a weaker version of the statement of Conjecture \ref{conj1}.  

We write $\data$ for the set of all conjugacy classes of pairs $(M,\delta)$ 
consisting of a Levi subgroup $M$ of $\G(\R)^1$ and a discrete series 
representation $\delta$ of $M^1$, where $M=A_M\times M^1$ and
$A_M$ is the largest central subgroup of $M$ isomorphic to a power of 
$\R^{>0}$. For any $\underline{\delta} \in \data$ let
$\Pi(\G(\R)^1)_{\underline{\delta}}$ be the set of all irreducible unitary
representations which arise by the Langlands quotient construction
from the irreducible constituents of $I_{M}^{L}(\delta)$ for Levi subgroups 
$L \supset M$. Here, $I_{M}^L$ denotes (unitary) induction from an arbitrary 
parabolic subgroup of $L$ with Levi subgroup $M$ to $L$.

\begin{definition}
Let ${\mathfrak M}$ be a set of Borel measures on $\Pi(\G(\R)^1)$. We call 
${\mathfrak M}$ {\it polynomially bounded} (PB), if for all 
$\underline\delta\in{\mathcal D}$ there exist $N_{\underline\delta}>0$ such that
\[
\mu\left(\{\pi\in\Pi(\G(\R)^1)_{\underline\delta}\colon|\lambda_\pi|\le R\}\right)
\ll_{\underline\delta}(1+R)^{N_{\underline\delta}}
\]
for all $\mu\in{\mathfrak M}$ and $R>0$. 
\end{definition}

Now consider the measures $\mu_K$ defined by \eqref{adelic-meas}. Let 
$\M\in \cL$ and denote by $\K_M(N)$ the congruence subgroups of $\M(\A_{\fin})$,
defined by \eqref{congruence}. Denote by $\mu^{\M}_{\K_{M}(N)}$ the measure 
defined by \eqref{measure4} with $\M$ in place of $\G$.
Then the key result is the following lemma.
\begin{lem}
Suppose that $\G$ satisfies (TWN) and (BD). Then for each $\M\in\cL$, the
collection of measures $\{\mu^{\M}_{\K_{M}(N)}\}$, $N\in\N$, is polynomially
bounded.
\end{lem}
This has the consequence that if $\G$ satisfies (TWN) and (BD), then for 
every $\M\neq\G$
and $f\in C^\infty_{c,\fin}(\G(\R)^1)$ we have
\[
J_{\spec,M}(f\otimes {\bf 1}_{\K(N)})\to 0
\]
as $N\to\infty$. Thus by Theorem \ref{thm-specexpand} it follows that if 
$\G$ satisfies (TWN) and (BD), then for every $f\in C^\infty_{c,\fin}(\G(\R)^1)$
we have
\[
J_{\spec}(f\otimes {\bf 1}_{\K(N)})-\tr R_{\di}(f\otimes {\bf 1}_{\K(N)})\to 0
\]
for $n\to\infty$. Thus  the spectral limit property is satisfied in this case.
By Theorems \ref{mainglobal} and \ref{thm-bd}, the groups $\GL(n)$ and 
$\SL(n)$ satisfy (TWN) and (BD) and therefore, the spectral limit property
holds for $\GL(n)$ and $\SL(n)$.

To deal with the geometric limit property we use the coarse geometric 
expansion 
\begin{equation} 
J^T (h) = \sum_{\mathfrak{o} \in \mathcal{O}} J^T_\mathfrak{o} (h),\quad h\in 
C_c^\infty(\G(\A)^1),
\end{equation}
(see \eqref{geometricside} for the notation). Write 
$J_\mathfrak{o}(f)=J_\mathfrak{o}^{T_0}(f)$, which depends only on $\M_0$ and $\K$. 
Let $J^T_{\unip}$ be the contribution of the unipotent elements of $G(\Q)$ to 
the trace formula \eqref{geometricside}, which is a polynomial in 
$T \in \aaa_{M_0}$ of degree at most $d_0 = \dim \aaa^G_{M_0}$ \cite{Ar7}.
It can be split into the contributions of the finitely many 
$G(\bar \Q)$-conjugacy classes of unipotent elements of $\G(\Q)$. It is well 
known ([ibid., Corollary 4.4]) that the contribution of the unit element is 
simply the constant polynomial $\vol (\G (\Q) \bs \G (\A)^1) h (1)$.
Write
\[
J^T_{\unip - \{ 1 \}} (h) = J^T_{\unip} (h) - \vol (\G (\Q) \bs \G (\A)^1) h (1), 
\quad h \in C^\infty_c (\G (\A)^1).
\]
Define the distributions $J_{\unip}$ and $J_{\unip - \{ 1 \}}$ as 
$J_{\unip}^{T_0}$ and $J_{\unip - \{ 1 \}}^{T_0}$, respectively. Since the groups 
$\K(N)$ form a neighborhood basis of the 
identity element in $\G(s\A_{\fin})$, it is easy to see that for a given
$h\in C_c^\infty(\G(\A)^1)$, for all but finitely many $N$ one has
\begin{equation}\label{unipotent1}
J(h\otimes {\bf 1}_{\K(N)})=J_{\unip}(h\otimes {\bf 1}_{\K(N)}).
\end{equation}

For any compact subset $\Omega\subset G(\R)^1$ we write 
$C^\infty_\Omega(\G(\R)^1)$ for the Fr\'echet space of
all smooth functions on $G(\R)^1$ supported in $\Omega$ equipped with the 
seminorms $\sup_{x \in \Omega} \abs{(Xh)(x)}$,
where $X$ ranges over the left-invariant differential operators on 
$\G(\R)$.
The key result is the following proposition.
\begin{prop} \label{unipproposition}
For any compact subset $\Omega\subset \G(\R)^1$ there exists a seminorm 
$||\cdot||$ on $C^\infty_\Omega(\G(\R)^1)$ such that
\[
\abs{J_{\unip - \{ 1 \}} (h \otimes {\bf 1}_{\K(N)})} \le 
\frac{(1 + \log(N))}{N} \| h\|
\]
for all $h \in C^\infty_\Omega (\G (\R)^1)$ and all $N\in\N$.
\end{prop}
The proof of Proposition \ref{unipproposition}
consists of a slight extension of Arthur's arguments in \cite{Ar7}.
Combining \eqref{unipotent1} and Proposition \ref{unipproposition} the
geometric limit property follows. This completes the proof of Theorem \ref{lm1}
for $F=\Q$. The case of a general $F$ is proved similarly. For details see
\cite{FLM2}.

\section{Analytic torsion and torsion in the cohomology of arithmetic groups }
\setcounter{equation}{0}

The theorem of DeGeorge and Wallach  on limit multiplicities for 
discrete series \cite{DW1} implies the statement 
\eqref{betti} on the approximation of $L^2$-Betti numbers by normalized Betti
numbers of finite covers \cite{AB2}. For towers of normal subgroups of finite 
index, L\"uck \cite{Lu1} proved this in the more general context of finite 
CW complexes. This is part of his study of the approximation of $L^2$-invariants
by their classical counterparts \cite{Lu2}. A more sophisticated spectral 
invariant is the analytic torsion introduced by Ray and Singer \cite{RS}. The 
study of the corresponding approximation
problem has interesting applications to the torsion in the cohomology of 
arithmetic groups. 
 
\subsection{Analytic torsion and $L^2$-torsion}
Let $X$ be a compact Riemannian manifold of dimension $n$  and  let
$\rho\colon \pi_1(X)\to\GL(V)$
a finite dimensional representation of its fundamental group. 
Let $E_\rho\to X$ be the flat vector bundle associated with $\rho$. Choose a
Hermitian fiber metric in $E_\rho$. Let $\Delta_p(\rho)$ be the Laplace 
operator on
$E_\rho$-valued $p$-forms with respect to the  metrics on $X$ and in $E_\rho$. 
It is
an elliptic differential operator, which is formally self-adjoint and 
non-negative.
Since $X$ is compact, $\Delta_p(\rho)$ has a pure discrete spectrum consisting
of sequence of eigenvalues $0\le\lambda_0\le\lambda_1\le\cdots\to\infty$ of 
finite multiplicity. Let
\begin{equation}
\zeta_p(s;\rho):=\sum_{\lambda_j>0}\lambda_j^{-s}
\end{equation}
be the zeta function of $\Delta_p(\rho)$. The series converges absolutely and
uniformly on compact subsets of the half-plane $\Re(s)>n/2$ and admits a 
meromorphic extension to $s\in\C$, which is holomorphic at $s=0$. Then the
Ray-Singer analytic torsion $T_X(\rho)\in\R^+$ is defined by
\begin{equation}\label{anal-tors}
T_X(\rho):=\exp\left(\frac{1}{2}\sum_{p=1}^n (-1)^p p\frac{d}{ds}\zeta_p(s;\rho)
\big|_{s=0}\right).
\end{equation}
It depends on the metrics on $X$ and $E_\rho$. However, if $\dim(X)$ is odd and
$\rho$ acyclic, which means that $H^*(X,E_\rho)=0$, then $T_X(\rho)$ is 
independent of the metrics \cite{Mu3}. The analytic torsion has a topological 
counterpart. This is
the Reidemeister torsion $T^{\topo}_X(\rho)$ (usually it is denoted by
$\tau_X(\rho)$), which is defined in terms of a 
smooth triangulation of $X$ \cite{RS}, \cite{Mu1}. It is known that for 
unimodular representations $\rho$
(meaning that $|\det\rho(\gamma)|=1$ for all $\gamma\in\pi_1(X)$) one has
the equality 
\begin{equation}\label{antor-reidem}
T_X(\rho)=T^{\topo}_X(\rho)
\end{equation}
\cite{Ch}, \cite{Mu1}. In the general case of a non-unimodular representation
the equality does not hold, but the defect can be described \cite{BMZ}.

Let $X_i\to X$, $i\in\N$, be sequence of finite 
coverings of $X$. Let $\inf(X_j)$ denote the injectivity radius of $X_j$ and
assume that $\inj(X_j)\to\infty$ as $j\to\infty$. Then the question is: Does
\begin{equation}\label{normanaltor}
\frac{\log T_{X_j}(\rho)}{\vol(X_j)}
\end{equation}
converge as $j\to\infty$ and if so, what is the limit? For a tower of normal
coverings and the trivial representation $\rho_0$ a conjecture of L\"uck  
\cite[Conjecture 7.4]{Lu2} states that the sequence \eqref{normanaltor} 
converges and the
limit is the $L^2$-torsion, first introduced by Lott \cite{Lo} and Mathei 
\cite{MV}. The $L^2$-torsion is defined as follows. Recall that the zeta
function $\zeta_p(s)$ can be expressed in terms of the heat operator
\[
\zeta_p(s)=\frac{1}{{\bf \Gamma}(s)}\int_0^\infty
(\Tr\left(e^{-t\Delta_p}\right)-b_p)t^{s-1}\;dt,
\]
where $b_p$ is the $p$-th Betti number and $\Re(s)>n/2$. Let 
$e^{-t\widetilde\Delta_p}$ be the heat operator of the Laplace operator
$\widetilde\Delta_p$ on $p$-forms on the universal covering 
$\widetilde X$ of $X$.
Let $\widetilde K_p(t,x,y)$ be the kernel of $e^{-t\widetilde\Delta_p}$.
Note that $\widetilde K_p(t,x,y)$ is a homomorphism of $\Lambda^p T^\ast_y(X)$
to $\Lambda^p T^\ast_x(X)$.
Let $F\subset \widetilde X$ be a fundamental domain for the action of 
$\Gamma:=\pi_1(X)$ on $\widetilde X$. Then the $\Gamma$-trace of 
$e^{-t\widetilde\Delta_p(\rho)}$ is defined as
\begin{equation}\label{gamma-trace}
\Tr_\Gamma\left(e^{-t\widetilde\Delta_p}\right):=\int_F \tr\widetilde K_p(t,x,x)\;
dx.
\end{equation} 
The $L^2$-Betti number $b_p^{(2)}$ is defined as
\[
b_p^{(2)}:=\lim_{t\to\infty}\Tr_\Gamma\left(e^{-t\widetilde\Delta_p}\right).
\]
In order to be able to define the Mellin transform of the $\Gamma$-trace one 
needs to know
the asymptotic behavior of $\Tr_\Gamma(e^{-t\widetilde\Delta_p})$ as 
$t\to 0$ and $t\to\infty$. 
Using a parametrix for the heat kernel which is pulled back from a parametrix
on $X$, one can show that for $t\to0$, $\Tr_\Gamma(e^{-t\widetilde\Delta_p})
$ has an asymptotic expansion similar to the compact case \cite{Lo}. For the 
large time behavior we need to introduce the Novikov-Shubin invariants
\begin{equation}\label{nov-shub}
\widetilde\alpha_p=\sup\left\{\beta_p\in [0,\infty)\colon 
\Tr_\Gamma\left(e^{-t\widetilde\Delta_p}\right)-b_p^{(2)}=O(t^{-\beta_p/2})
\;\;\text{as}\:\;t\to\infty\right\}
\end{equation}
Assume that $\widetilde\alpha_p>0$ for all $p=1,\dots,n$. Then the 
$L^{2}$- torsion
$T^{(2)}_X\in\R^+$ can be defined by
\begin{equation}\label{l2-torsion}
\begin{split}
\log T^{(2)}_X=\frac{1}{2}\sum_{p=1}^n (-1)^pp \biggl[
\frac{d}{ds}\biggl(\frac{1}{\Gamma(s)}\int_0^1
\Tr_\Gamma\biggl(&e^{-t\widetilde\Delta_p^\prime}\biggr)t^{s-1}\,dt\biggr)\bigg|_{s=0}\\
&+\int_1^\infty t^{-1}\Tr_\Gamma\left(e^{-t\widetilde\Delta_p^\prime}\right)\;dt\biggr],
\end{split}
\end{equation}
where $\widetilde\Delta_p^\prime$ denotes the restriction of 
$\widetilde\Delta_p$ to the orthogonal complement of 
$\ker \widetilde\Delta_p$ and
the first integral is defined near $s=0$ by analytic continuation.
This definition can be generalized to all finite dimensional representations
$\rho$ of $\Gamma$, if the corresponding Novikov-Shubin invariants are all
positive. Then the $L^2$-torsion $T^{(2)}_X(\rho)$ is defined as in 
\eqref{l2-torsion}. If there exists $c>0$ such that the spectrum
of $\Delta_p(\rho)$ is bounded from below by $c$, then the integral
\[
\int_0^\infty \Tr_\Gamma\left(e^{-t\widetilde\Delta_p(\rho)}\right)t^{s-1}\;dt
\]
converges for $\Re(s)>n/2$ and admits a meromorphic continuation to $\C$
which is holomorphic at $s=0$. Thus, if there is a positive lower bound 
of the spectrum of all $\Delta_p(\rho)$, $p=1,\dots,n$, then $T^{(2)}_X(\rho)$
can be defined in the usual way by
\[
\log T^{(2)}_X(\rho)=\frac{1}{2}\sum_{p=1}^n (-1)^p p\frac{d}{ds}
\left(\frac{1}{\Gamma(s)}
\int_0^\infty\Tr_\Gamma\left(e^{-t\widetilde\Delta_p(\rho)}\right)t^{s-1}\,dt\right)
\bigg|_{s=0}.
\]
Let $\Gamma=\pi_1(X,x_0)$ and let $(\Gamma_i)_{i\in\N_0}$ be a tower of normal
subgroups of finite index of $\Gamma=\Gamma_0$. Let 
$X_i=\Gamma_i\bs \widetilde X$, $i\in\N_0$, be the corresponding covering of
$X$. Let $T_X$ and $T^{(2)}_X$ denote the analytic torsion and $L^2$-torsion
with respect to the trivial representation.
W. L\"uck \cite[Conjecture 7.4]{Lu2} has made the following conjecture.
\begin{conjecture} For every closed Riemannian manifold $X$ the $L^2$-torsion
$T^{(2)}_X$ exists and for a sequence of coverings $(X_i\to X)_{i\in\N}$ as above 
one has
\[
\lim_{i\to\infty}\frac{\log T_{X_i}}{[\Gamma\colon\Gamma_i]}=
\log T^{(2)}_{X}.
\]
\end{conjecture}
One is tempted to make this conjecture for any finite dimensional representation
$\rho$. 

\subsection{Compact locally symmetric spaces}

Now we turn to the locally symmetric case. Let $X=\Gamma\bs \widetilde X$,
where $\widetilde X=G/K$ is a Riemannian symmetric space of non-positive 
curvature and $\Gamma\subset G$ is a discrete, torsion free, cocompact subgroup.
Let $\tau$ be an irreducible finite dimensional complex representation of $\G$.
Let $E_\tau\to X$ be the flat vector bundle associated to the representation
$\tau|_\Gamma$ of $\Gamma$. By \cite{MM}, $E_\tau$ can be equipped with a 
canonical Hermitian fiber metric, called admissible, which is unique up to 
scaling. Let $\Delta_p(\tau)$ be the Laplace operator on $p$-forms with values
in $E_\tau$, with respect to the choice of any admissible fiber metric in
$E_\tau$.  Let $T_X(\tau)$ be the corresponding analytic torsion. 
Let $\widetilde\Delta_p(\tau)$ be the
Laplace operator on $\widetilde E_\tau$-valued $p$-forms on $\widetilde X$.
Let $\widetilde E_\tau\to\widetilde X$ be the homogeneous vector bundle 
defined by $\tau|_K$. By \cite{MM} there is a canonical isomorphism
\[
E_\tau\cong \Gamma\bs\widetilde E_\tau
\]
and the metric on $E_\tau$ is induced by the homogeneous metric on $\widetilde 
E_\tau$. Thus
\begin{equation}\label{sections}
C^\infty(\widetilde X,\widetilde E_\tau)\cong (C^\infty(G)\otimes V_\tau)^K.
\end{equation}
Let $R$ be the right regular representation of $G$ in $C^\infty(G)$ and let 
let $R(\Omega)$ be the operator in $(C^\infty(G)\otimes V_\tau)^K$ induced by
the Casimir element. Then with respect to the isomorphism \eqref{sections} we
have
\[
\widetilde\Delta_p(\tau)=-R(\Omega)+\lambda_\tau\Id
\]
(see \cite{MM}).  This implies that the heat operator
$e^{-t\widetilde\Delta_p(\tau)}$ is a convolution operator given by a kernel
\[
H_t^{p,\tau}\colon G\to \End(\Lambda^p\pf^\ast\otimes V_\tau).
\]
Let $h_t^{p,\tau}\in C^\infty(G)$ be defined by $h_t^{p,\tau}(g)=\tr H_t^{p,\tau}(g)$,
$g\in G$. Then it follows from \eqref{gamma-trace} that
\begin{equation}\label{l2-trace}
\Tr_\Gamma\left(e^{-t\widetilde\Delta_p(\tau)}\right)=\vol(X)h_t^{p,\tau}(1).
\end{equation}
Now one can use the Plancherel theorem to compute $h_t^{p,\tau}(1)$ and
determine its asymptotic bahavior as $t\to 0$ and $t\to\infty$. For the
trivial representation this was carried out in \cite{Ol} and for strongly 
acyclic $\tau$ in \cite{BV}. So let $\widetilde\Delta_p(\tau)^\prime$ be the 
restriction of $\widetilde\Delta_p(\tau)$ to the orthogonal complement of the
kernel of $\widetilde\Delta_p(\tau)$. Now let
\begin{equation}
\widetilde\alpha_p(X,\tau):=\sup\left\{\beta_p\in [0,\infty)\colon 
\Tr_\Gamma\left(e^{-t\widetilde\Delta_p(\tau)^\prime}\right)=O(t^{-\beta_p/2})
\;\;\text{as}\:\;t\to\infty\right\}, 
\end{equation}
$p=0,\dots,n$, be the twisted Novikov-Shubin invariants. 
Assume that $\widetilde\alpha_p(X,\tau)>0$, $p=0,\dots,n$. Then the 
$L^2$-torsion $T^{(2)}_X(\tau)$ is defined. By \cite[Theorem 1.1]{Ol} this is 
the case for the trivial representation. Furthermore, if $\tau$ is strongly 
acyclic, then $\widetilde\alpha_p(X,\tau)=\infty$ for all $p$. 
Using the definition of the $L^2$-torsion, it follows that
\begin{equation}\label{l2-tor2}
\log T_X^{(2)}(\tau)=\vol(X)t^{(2)}_{\widetilde X}(\tau),
\end{equation}
where $t^{(2)}_{\widetilde X}(\tau)$ is a constant that depends only on 
$\widetilde X$ and $\tau$. 

Now let $(\Gamma_j)$ be sequence of torsion free cocompact lattices in $G$.
Let $X_j=\Gamma_j\bs\widetilde X$ and assume that $\inj(X_j)\to \infty$ if 
$j\to\infty$. A representation $\tau\colon G\to\GL(V)$ is called {\it strongly
acyclic}, if there is $c>0$ such that the spectrum of $\Delta_{X_j,p}(\tau)$ is
contained in $[c,\infty)$ for all $j\in\N$ and $p=0,\dots,n$. 

Now let $\G$ be a connected semisimple algebraic $\Q$-group. Let $G=\G(\R)$.  
Then it is proved in \cite{BV} that strongly acylic representations exist.
For such representations Bergeron and Venkatesh 
\cite[Theorem 4.5]{BV} established the following theorem.
\begin{theo}\label{thm-torconv}
Let $\tau\colon G\to\GL(V)$ be strongly acyclic. Then
\begin{equation}\label{converg1}
\lim_{j\to\infty}\frac{\log(T_{X_j}(\tau))}{\vol(X_j)}=t^{(2)}_X(\tau),
\end{equation}
where $X_j=\Gamma_j\bs \widetilde X$ and $\inj(X_j)\to\infty$ as $j\to\infty$.
\end{theo}
The number $t^{(2)}_X(\tau)$ can be computed using the Plancherel theorem. Let
$\delta(G)=\rk(G)-\rk(K)$ be the fundamental rank or ``deficiency'' of $G$.
By \cite[Proposition 5.2]{BV} one has 
\begin{prop}\label{l2-tors3}
If $\delta(G)\neq 1$, then $t^{(2)}_X(\tau)=0$. For $\delta(G)=1$ one has
\[
(-1)^{\frac{\dim\widetilde X-1}{2}} t^{(2)}_X(\tau)>0.
\]
\end{prop}
We note that the simple Lie groups $G$ with $\delta(G)=1$ are $\SL_3(\R)$ and 
$\SO(p,q)$ with $pq$ odd, especially $G=\SO^0(2m+1,1)$ is a group with
fundamental rank 1. 
 
Next we briefly recall the main steps of the proof of Theorem \ref{thm-torconv}.
To indicate the dependence 
of the heat operator and other quantities on the covering $X_j$, we use the 
subscript $X_j$. The uniform spectral gap at $0$ implies that there exist
constants $C,c>0$ such that for all $p=0,\ldots,n$, $j\in\N$ and $t\ge 1$ one 
has
\begin{equation}\label{est-trace2}
\Tr\left(e^{-t\Delta_{X_j,p}(\tau)}\right)\le C e^{-tc}\vol(X_j)
\end{equation}
(see \cite{BV}). This is the key result that makes the method to work. Let
\begin{equation}\label{analtor1}
K_{X_j}(t,\tau):=\frac{1}{2}\sum_{p=1}^n(-1)^pp
\Tr\left(e^{-t\Delta_{X_j,p}(\tau)}\right).
\end{equation}
Using \eqref{est-trace2} it follows that the analytic torsion can be defined
by
\begin{equation}\label{analtor2}
\log T_{X_j}(\tau)=\frac{d}{ds}\left(\frac{1}{\Gamma(s)}\int_0^\infty
K_{X_j}(t,\tau)t^{s-1}\;dt\right)\bigg|_{s=0}.
\end{equation}
Let $T>0$. Then we can split the integral and rewrite the right hand side as
\[
\log T_{X_j}(\tau)=\frac{d}{ds}\left(\frac{1}{\Gamma(s)}\int_0^T
K_{X_j}(t,\tau)t^{s-1}\;dt\right)\bigg|_{s=0}+\int_T^\infty K_{X_j}(t,\tau)t^{-1}\;dt.
\]
By \eqref{est-trace2} there exist $C,c>0$ such that
\begin{equation}\label{part-mellin}
\frac{1}{\vol(X_j)}\left|\int_T^\infty K_{X_j}(t,\tau)t^{-1}\;dt\right|\le C
e^{-cT}
\end{equation}
for all $j\in\N_0$ and $T>1$. To deal with the first term one can use the 
Selberg trace formula. Put
\[
k_t^\tau:=\frac{1}{2}\sum_{p=1}^n(-1)^pp h_t^{p,\tau}.
\]
Then the Selberg trace formula gives
\[
K_{X_j}(t,\tau)=\vol(X_j)k_t^\tau(1)+H_{X_j}(k_t^\tau),
\]
where $H_{X_j}(k_t^\tau)$ is the contribution of the hyperbolic conjugacy classes.
Using \eqref{l2-trace} and the definition of $k_t^\tau$, it follows that
\[
\frac{d}{ds}\left(\frac{1}{\Gamma(s)}\int_0^T k_t^\tau(1)t^{s-1}\;dt\right)
\bigg|_{s=0}=t^{(2)}_{\widetilde X}(\tau)+O\left(e^{-cT}\right)
\]
as $T\to\infty$. Regrouping the terms of the hyperbolic contribution
$H_{X_j}(k_t^\tau)$ as in \eqref{corwin2} it follows that the corresponding 
integral divided by $\vol(X_j)$ converges to 0 as $j\to\infty$. This 
proves the theorem. 

One expects Theorem \ref{thm-torconv} to be true in general.  
However, if there is no spectral 
gap at zero, one cannot argue as above. The key problem is to control the
small eigenvalues as $j\to\infty$. Sufficient conditions 
on the behavior of the small eigenvalues are discussed in \cite{Lu2}
and in the 3-dimensional case also in \cite{BSV}. 

In view of the potential applications to the cohomology of arithmetic groups,
discussed in the next section, it is very desirable to extend Theorem 
\ref{thm-torconv} to the non-compact case. The first problem one faces is that 
the corresponding Laplace operators have a nonempty continuous spectrum and
therefore, the heat operators are not trace class and the analytic torsion can
not be defined as above. This problem has been studied by Raimbault \cite{Ra1}
for hyperbolic 3-manifolds and in \cite{MP2} for hyperbolic manifolds of any 
dimension. 

So let $G=\SO^0(n,1)$, $K=\SO(n)$ and $\widetilde X=G/K$. Equipped with a
suitably normalized $G$-invariant metric, $\widetilde X$ becomes isometric to
the $n$-dimensional hyperbolic space $\bH^n$.  Let $\Gamma\subset G$ be a 
torsion free lattice. Then $X=\Gamma\bs\widetilde X$ is an oriented 
$n$-dimensional hyperbolic manifold of finite volume.
As above, let $\tau\colon G\to\GL(V)$ be a finite dimensional complex
representation of $G$. The first step is to define a regularized trace of the 
heat operators $e^{-t\Delta_p(\tau)}$. To this end one uses an appropriate height
function to truncate $X$ at sufficient high level $Y>Y_0$ to get a compact
manifold $X(Y)\subset X$ with boundary $\partial X(Y)$, which consists of
a disjoint union of $n-1$-dimensional tori. Let $K^{p,\tau}(t,x,y)$ be the kernel
of the heat operator $e^{-t\Delta_p(\tau)}$. Using the spectral resolution of
$\Delta_p(\tau)$, it follows that there exist $\alpha(t)\in\R$ such that
$\int_{X(Y)}\tr K^{p,\tau}(t,x,x)\;dx-\alpha(t)\log Y$ has a limit as $Y\to\infty$.
Then we define the regularized trace as
\begin{equation}\label{reg-trace}
\Tr_{\reg}\left(e^{-t\Delta_p(\tau)}\right):=\lim_{Y\to\infty}
\left(\int_{X(Y)}\tr K^{p,\tau}(t,x,x)\;dx-\alpha(t)\log Y\right).
\end{equation}
We note that the regularized trace is not uniquely defined. It depends on the 
choice of truncation parameters on the manifold $X$. However, if 
$X_0=\Gamma_0\bs\bH^n$ is given and if truncation parameters on $X_0$ are fixed,
then every finite covering $X$ of $X_0$ is canonically equipped with 
truncation parameters, namely one simply pulls back the height function on $X_0$
to a height function on $X$ via the covering map. 

Let $\theta$ be the Cartan involution of $G$ with respect to $K=\SO(n)$. Let
$\tau_\theta=\tau\circ\theta$. If $\tau\not\cong\tau_\theta$, it can be shown
that $\Tr_{\reg}\left(e^{-t\Delta_p(\tau)}\right)$ is exponentially decreasing as
$t\to\infty$ and admits an asymptotic expansion as $t\to 0$. Therefore, the
regularized zeta function $\zeta_{\reg,p}(s;\tau)$ of $\Delta_p(\tau)$ can be 
defined
as in the compact case by
\begin{equation}
\zeta_{\reg,p}(s;\tau):=\frac{1}{\Gamma(s)}\int_0^\infty\Tr_{\reg}
\left(e^{-t\Delta_p(\tau)}\right)t^{s-1}\;dt.
\end{equation}
The integral converges absolutely and uniformly on compact subsets of the 
half-plane $\Re(s)>n/2$ and admits a meromorphic extension to the whole complex 
plane, which is holomorphic at $s=0$. So in analogy with the compact case, the 
regularized analytic torsion $T_X(\tau)\in\R^+$ can be defined by the 
same formula \eqref{anal-tors}. 

In even dimension the analytic torsion is rather trivial. Therefore, we assume
that $n=2m+1$. Furthermore, for  technical reasons we assume that 
every lattice $\Gamma
\subset G$ satisfies the following condition: For every $\Gamma$-cuspidal
parabolic subgroup $P$ of $G$ one has
\begin{equation}\label{intersect}
\Gamma\cap P=\Gamma\cap N_P,
\end{equation}
where $N_P$ denotes the unipotent radical of $P$.  
Let $\Gamma_0$ be a fixed lattice in $G$ and
let $X_0=\Gamma_0\bs \widetilde X$. Let $\Gamma_j$, $j\in\N$, be a sequence 
of finite index torsion free subgroups of $\Gamma_0$. This sequence is called
to be {\it cusp uniform}, if the tori which arise as cross sections
of the cusps of the manifolds $X_J:=\Gamma_j\bs \widetilde X$ satisfy some
uniformity condition (see \cite[Definition 8.2]{MP2}). 

The following theorem and its corollaries are established in \cite{MP2}.
One of the main results of \cite{MP2} is the following theorem 
which  
may be regarded as an analog of  Theorem \ref{thm-torconv} for oriented finite 
volume hyperbolic manifolds. 
\begin{theo}\label{thm-finitevol}
Let $\Gamma_0$ be a lattice in $G$ and let $\Gamma_i$, $i\in\mathbb{N}$, be a
sequence 
of finite-index normal subgroups which is cusp uniform and such that 
each $\Gamma_i$, $i\geq 1$, is torsion-free and satisfies \eqref{intersect}. If
$\lim_{i\to\infty}[\Gamma_0:\Gamma_i]=\infty$ and if each
$\gamma_0\in\Gamma_0-\{1\}$ 
only belongs to finitely many $\Gamma_i$, then for each $\tau$ with
$\tau\neq\tau_\theta$ one has  
\begin{align}\label{eqtheo}
\lim_{i\to\infty}\frac{\log T_{X_i}(\tau)}{[\Gamma:\Gamma_i]}=
t^{(2)}_{\bH^n}(\tau)\vol(X_0).
\end{align}
In particular, if under the same assumptions $\Gamma_i$ is a tower of normal
subgroups,
i.e. $\Gamma_{i+1}\subset
\Gamma_i$ for each $i$ and 
$\cap_i\Gamma_i=\{1\}$, then \eqref{eqtheo} holds.
\end{theo}

For hyperbolic 3-manifolds, Theorem \ref{thm-finitevol} was proved by J. 
Raimbault
\cite{Ra1} under additional assumptions on the intertwining operators. We
emphasize that the above theorem holds without any additional assumptions.

Now we specialize to arithmetic groups. First consider 
$\Gamma_0:=\SO^0(n,1)(\Z)$. Then $\Gamma_0$ is a lattice in $\SO^0(n,1)$. For 
$q\in\N$ let $\Gamma(q)$ be the principal congruence subgroup of $\Gamma_0$
of level $q$. Using a result of Deitmar and Hoffmann \cite{DH}, it follows
that the family of principal congruence subgroups $\Gamma(q)$ is cusp 
uniform \cite[Lemma 10.1]{MP2}. Thus Theorem \ref{thm-finitevol} implies the
following corollary (see \cite[Corollary 1.3]{MP2}).
\begin{corollary}\label{congruence1} For any finite dimensional irreducible 
representation $\tau$ of $\SO^0(n,1)$
with $\tau\not\cong\tau_\theta$ the principal congruence subgroups $\Gamma(q)$,
$q\ge 3$, of $\Gamma_0:=\SO^0(n,1)(\Z)$ satisfy
\[
\lim_{q\to\infty}\frac{\log T_{X_q}(\tau)}{[\Gamma\colon\Gamma(q)]}=
t^{(2)}_{\bH^n}(\tau) \vol(X_0),
\]
where $X_q:=\Gamma(q)\bs \bH^n$ and $X_0:=\Gamma_0\bs\bH^n$.
\end{corollary}
We recall that by Proposition \ref{l2-tors3} we have 
$(-1)^{\frac{n-1}{2}}t^{(2)}_{\bH^n}(\tau)>0$.

Next we consider  the 3-dimensional case. We note that every lattice 
$\Gamma\subset\SO^0(3,1)$ can be lifted to a lattice $\Gamma^\prime\subset 
\Spin(3,1)$. Moreover, recall that there is a natural isomorphism 
$\Spin(3,1)\cong\SL_2(\C)$. If $\rho$ is the standard representation of 
$\SL_2(\C)$ on $\C^2$, then the finite dimensional irreducible representations
of $\SL_2(\C)$ are given by $\Sym^p\rho\otimes\Sym^q\bar\rho$, $p,q\in\N$, where
$\Sym^k$ denotes the $k$-th symmetric power and $\bar\rho$ denotes the
complex conjugate representation to $\rho$. One has $(\Sym^p\rho\otimes
\Sym^q\bar\rho)_\theta=\Sym^q\rho\otimes\Sym^p\bar\rho$. For $D\in\N$ square free
let ${\mathcal O}_D$ be the ring of integers of the imaginary quadratic field
$\Q(\sqrt{-D})$ and let $\Gamma(D):=\SL_2({\mathcal O}_D)$. Then $\Gamma(D)$ is 
a lattice in $\SL_2(\C)$. If $\af$ is a non-zero ideal in ${\mathcal O}_D$, let
$\Gamma(\af)$ be the associated principal congruence subgroup of level $\af$.
Then Theorem \ref{thm-torconv} implies the following corollary (see 
\cite[Corollary 1.4]{MP2}.

\begin{corollary}\label{bianchigr1}
Let $D\in\N$ be square free. 
Let $\af_i$ be a sequence of non-zero ideals in ${\mathcal O}_D$ such that each
$N(\af_i)$ is sufficiently large and such that $\lim_{i\to\infty}N(\af_i)=\infty$.
Put  $X_D:=\Gamma(D)\bs \bH^3$ and $X_i:=\Gamma(\af_i)\bs\bH^3$. Let
 $\tau=\Sym^p\rho\otimes\Sym^q\bar\rho$ with $p\neq q$. Then one has
\[
\lim_{i\to\infty}\frac{\log T_{X_i}(\tau)}{[\Gamma(D)\colon \Gamma(\af_i)]}=
t^{(2)}_{\bH^3}(\tau)\vol(X_D).
\]
\end{corollary}  

\subsection{Applications to the cohomology of arithmetic groups - 
the cocompact case}

Theorem \ref{thm-torconv} has interesting consequences for the cohomology of 
arithmetic groups.  Let $\Gamma\subset G$ be a discrete, torsion free, 
cocompact  subgroup. Let $\tau\colon G\to \GL(V)$ be a finite dimensional real
representation and let $E\to X$ be the associated vector bundle.
Choose a fiber metric $h$ in $E$. 
 Assume that there exist a $\Gamma$-invariant lattice 
$M\subset V$. Let $\cM$ be the associated local system of free $\Z$-modules 
over $X$. Then we have $E=\cM\otimes \R$. Let $H^\ast(X,\cM)$ be the
cohomology of $X$ with coefficients in $\cM$. Each $H^q(X,\cM)$ is a
finitely generated $\Z$-module. Let $H^q(X,\cM)_{\tors}$ be the torsion
subgroup and
\[
H^q(X;\cM)_{\free}=H^q(X,\cM)/H^q(X,\cM)_{\tors}.
\]
We identify  $H^q(X,\cM)_{\free}$ with a subgroup of $H^q(X,E)$.
Let $\langle\cdot,\cdot\rangle_q$ be the inner product in $H^q(X,E)$ 
induced by the $L^2$-metric on ${\mathcal H}^q(X,E)$. 
 Let $e_1,...,e_{r_q}$ be a basis of $H^q(X,\cM)_{\free}$ and let $G_q$ be the
Gram matrix with entries $\langle e_k,e_l\rangle$. Put
\[
R_q(\tau,h)=\sqrt{|\det G_q|}, \quad q=0,...,n.
\]
Define the ``regulator'' $R(\tau,h)$ by
\begin{equation}\label{regulator1}
R(\tau,h)=\prod_{q=0}^n R_q(\tau,h)^{(-1)^q}.
\end{equation}
Recall that the  Reidemeister torsion $T^{\topo}_X(\tau,h)$ depends on the metric
$h$ through the choice of an orthonormal basis in the 
cohomology $H^\ast(X,E)$, where the inner product in $H^\ast(X,E_\tau)$
is defined as above. The key result relating Reidemeister torsion and 
cohomology is the following proposition.
\begin{prop}\label{reidem1}
 Let $\tau$ be a unimodular representation of $\Gamma$ on a finite-dimensional
$\R$-vector space $V$. Let $M\subset V$ be a $\Gamma$-invariant lattice and let
$\cM$ be the associated local system of finitely generated free $\Z$-modules 
on $X$. Let $h$ be
a fiber metric in the flat vector bundle $E=\cM\otimes\R$. Then we have
\begin{equation}\label{reid-coho}
T^{\topo}_X(\tau,h)=R(\tau,h)\cdot \prod_{q=0}^n 
\left|H^q(X,\cM)_{\tors}\right|^{(-1)^{q+1}}.
\end{equation}
\end{prop}
Especially, if $\tau|_\Gamma$ is acyclic, i.e., if $H^\ast(X,E)=0$, then 
$T^{\topo}_X(\tau,h)$ is independent of $h$ and we denote it by $T^{\topo}_X(\tau)$.
Moreover, $R(\tau,h)=1$. Then $H^\ast(X,\cM)$ is a torsion group and one has
\[
T^{\topo}_X(\tau)=\prod_{q=0}^n \left|H^q(X,\cM)\right|^{(-1)^{q+1}}.
\]
Representations $\tau$ of $G$ which admit a $\Gamma$-invariant lattice arise 
in the following arithmetic situation. Let $\G$ be a semisimple algebraic group
defined over $\Q$ and let $G=\G(\R)$. Let $\Gamma\subset \G(\Q)$ be an
arithmetic subgroup. Let $V_0$ be 
a $\Q$-vector space and let $\rho\colon \G\to \GL(V_0)$ be a rational 
representation. Then there exists a lattice $M\subset V_0$ which is invariant
under $\Gamma$ and $V_0=M\otimes_\Z\Q$. Let $V=V_0\otimes_\Q\R$ and let
$\tau\colon G\to \GL(V)$ be the representation induced by $\rho$. Then 
$M\subset V$ is a $\Gamma$-invariant lattice. 

Assume that  $\Gamma\subset \G(\Q)$ is cocompact in $G$ 
(equivalently assume that $\G$
is anisotropic). Then it is proved in \cite{BV} that strongly acyclic arithmetic
$\Gamma$-modules $M$ exist. Assume that $\delta(G)=1$. Let $M$ be a strongly
acyclic arithmetic $\Gamma$-module. Then by 
\eqref{antor-reidem}, Theorem \ref{thm-torconv} and Proposition \ref{l2-tors3}
it follows that there exists a constant $C>0$, which depends on $G$ and $M$,
such that 
\begin{equation}
\lim_{j\to\infty}\sum_{k=0}^n(-1)^{k+\frac{\dim(\widetilde X)-1}{2}}
\frac{\log|H_k(X_j,\cM)|}{[\Gamma\colon\Gamma_j]}=C\vol(\Gamma\bs \widetilde X)
\end{equation}
(see \cite[(1.4.2)]{BV}). This implies the following theorem of Bergeron and 
Venkatesh.
\cite[Theorem 1.4]{BV}.
\begin{theo}\label{hom-growth}
Suppose that $\delta(\widetilde X)=1$. Then strongly acyclic arithmetic
$\Gamma$-modules exist. For any such module $M$,
\[
\liminf_{j}\sum_{k\equiv a\hskip-6pt\pmod 2}
\frac{\log|H_k(X_j,\cM)|}{[\Gamma\colon\Gamma_j]}\ge C
\vol(\Gamma\bs \widetilde X),
\] 
where  $a=(\dim(\widetilde X)-1)/2$ and $C>0$ depends only on $G$ and $M$.
\end{theo}
In Theorem \ref{hom-growth}, one cannot in general isolate the degree which
produces torsion. 
A conjecture of Bergeron and Venkatesh \cite[Conjecture 1.3]{BV} claims the 
following.
\begin{conjecture}
The limit
\[
\lim_{j\to\infty}\frac{\log|H_k(X_j,\cM)_{\tors}|}{[\Gamma\colon\Gamma_j]}
\]
exists for each $k$ and is zero unless $\delta(G)=1$ and 
$k=\frac{\dim(\widetilde X)-1}{2}$. In that case, it is always positive and 
equal
to a positive constant $C_{G,M}$, which can be explicitly described, times 
$\vol(\Gamma\bs\widetilde X)$.
\end{conjecture}
An example, for which this conjecture can be verified is $G=\SL(2,\C)$. 

If the representation $\tau$ of $G$ is not acyclic, various difficulties
occur. First of all, the spectrum of the Laplace operators has no positive
lower bound which causes the problem with the small eigenvalues discussed
above in the context of analytic torsion. Secondly the regulator $R(\tau,h)$
is in general nontrivial. It turns out to be rather difficult to control the 
growth of the regulator. Of particular interest is the case of the trivial
representation, i.e., the integer homology $H_k(X_j,\Z)$.
The 3-dimensional case has been studied in \cite{BSV}. In this paper the
authors discuss  conditions which imply that the results of \cite{BV} on
strongly acyclic local systems can be extended to the case of the trivial 
local system. There are conditions on the cohomology and the spectrum of the
Laplace operator on 1-Forms. The conditions on the spectrum are as follows. 
Let $(\Gamma_i)_{i\in\N}$ be a sequence of cocompact congruence subgroups of a 
fixed arithmetic subgroup $\Gamma\subset \SL(2,\C)$. Let $X_i=\Gamma_i\bs
\bH^3$ and put $V_i:=\vol(X_i)$. Let $\lambda^{(i)}_j$ $j\in\N$, be the 
eigenvalues of the Laplace operator on $1$-forms of $X_i$. Assume:
\begin{enumerate}
\item[1)] For every $\varepsilon>0$ there exists $c>0$ such that
\[
\limsup_{i\to\infty}\frac{1}{V_i}\sum_{0<\lambda^{(i)}_j\le c} |\log\lambda^{(i)}
_j|\le \varepsilon.
\]
\item[2)] $b_1(X_i,\Q)=o(\frac{V_i}{\log V_i})$.
\end{enumerate}
Let $T_{X_i}$ be the analytic torsion with respect to the trivial local system.
As shown in \cite{BSV},   conditions 1) and 2) imply that
\[
\frac{\log T_{X_i}}{V_i}\longrightarrow t^{(2)}_{\bH^3}=-\frac{1}{6\pi}, \quad
i\to\infty.
\]
Unfortunately, it seems to be difficult to verify 1) and 2). The other problem 
is to estimate the growth of the regulator (see \cite{BSV}). We note that 
condition 1) is equivalent to the following condition $1^\prime$). 
\begin{enumerate}
\item[$1^\prime$)] Let $d\mu_{1}$ be the spectral measure of 
$\widetilde\Delta_1$. For every $c>0$ one has
\[
\frac{1}{V_i}\sum_{0<\lambda_j^{(i)}\le c}\log \lambda_j^{(i)}\longrightarrow 
\int_0^c\log\lambda\; d\mu_1(\lambda),\quad i\to\infty.
\]
\end{enumerate}
There is a certain similarity with the limit multiplicity problem.

Finally we note that there is related work by Calegari and Venkatesh 
\cite{CaV} who use 
analytic torsion to compare torsion in the cohomology of different arithmetic
subgroups of $\SL(2,\C)$ and establish a numerical form of a Jacquet-Langlands
correspondence in the torsion case.

\subsection {The finite volume case}

Many important arithmetic groups are not cocompact. So it is desirable to 
extend the results of the previous section to the finite volume case. 
In order to achieve this one has to deal with the following problems. 
\begin{enumerate}
\item[1)] Define an appropriate  regularized version $T_X^{\reg}(\rho)$ of the 
analytic torsion for a
finite volume locally symmetric space $X=\Gamma\bs\widetilde X$ 
and establish the analog of \eqref{converg1}. So  
let  $\Gamma_j\subset \Gamma$ be a sequence of subgroups of finite index and
$X_j:=\Gamma_j\bs \widetilde X$, $j\in\N$. Assume that  $\vol(V_j)\to\infty$. 
Under appropriate additional assumptions on the sequence $(\Gamma_j)_{j\in\N}$ 
one has to show that
\[
\lim_{j\to\infty}\frac{\log T_{X_j}^{\reg}(\rho)}{\vol(X_j)}=
t^{(2)}_{\widetilde X}(\rho).
\]
\item[2)] Show that $T_X^{\reg}(\rho)$ has a topological counterpart 
$T_X^{\topo}(\rho)$,
possibly the Reidemeister torsion of an intersection complex. 
\item[3)] If $E_\rho$ is arithmetic, i.e., if there is a local system
of finite rank free $\Z$-modules $\cM$ over $X$ such that $E_\rho=\cM\otimes\R$,
establish an analog of \eqref{reid-coho}.
\item[4)] Estimate the growth of the regulator.
\end{enumerate}
For hyperbolic manifolds 1) has been proved in \cite{Ra1} in the 
3-dimensional case and in
\cite{MP1} and \cite{MP2} in general. It would be very interesting to extend
these results to the higher rank case. $\SL(3,\R)$ seems to be doable.

Raimbault \cite{Ra2} has studied 2) in the 3-dimensional case  
and established a kind of asymptotic equality of analytic and Reidemeister
torsion, which is sufficient for the present purpose. Of course, the goal is to 
prove an exact equality. For hyperbolic manifolds there is some recent 
progress \cite{AR}. Unfortunately, this paper does not cover the relevant 
flat bundles. The method requires that the flat bundle can extended to the 
boundary at infinity. This is not the case for the flat bundles which arise
from representations of $G$ by restriction to $\Gamma$.
J. Pfaff \cite{Pf} has established a gluing formula for the
regularized analytic torsion of a hyperbolic manifold, which reduces the
problem to the case of a cusp.

4) has been studied by Raimbault \cite{Ra2} for 3-dimensional 
hyperbolic manifolds.  It turns out to be very difficult. The real cohomology
never vanishes. There is always the part of the cohomology coming from the
boundary. This is the Eisenstein cohomology introduced by Harder \cite{Ha}.
These cohomology classes are represented  by Eisenstein classes, which are 
rational cohomology classes. The problem is to
estimate the denominators of the Eisenstein classes which seems to be a hard
problem.

\end{document}